\documentclass[twoside,11pt]{amsart}
\usepackage{amsmath,amssymb}

\newcommand{\CP}{\mathbb{CP}}
\newcommand{\bCP}{\overline{\mathbb{CP}}}
\newcommand{\ta}{\tilde{a}}
\newcommand{\tb}{\tilde{b}}
\newcommand{\tx}{\tilde{x}}
\newcommand{\ty}{\tilde{y}}
\newcommand{\tF}{\tilde{F}}

\def \eu{{\text{e}}}
\def \sign{{\text{sign}}}
\newcommand{\CPb}{\overline{\mathbb{CP}}{}^{2}}

\newcommand{\BZ}{\mathbb{Z}}

\newcommand{\TryPackage}[3]{\IfFileExists{#1.sty}{\usepackage{#1}#2}{#3}
}
\TryPackage{mathrsfs}{\renewcommand{\mathcal}{\mathscr}}{%
        \TryPackage{eucal}{}{}}

\newcommand{\al}{\alpha}
\newcommand{\be}{\beta}

\newcommand{\ga}{\gamma}

\newcommand{\ZZ}{{\mathbb Z}}
\newcommand{\RR}{{\mathbb R}}

\usepackage{latexsym}
\usepackage{amsmath}
\usepackage{amsthm}
\usepackage{amscd}
\usepackage{xypic} %This is XY-pic.
\usepackage{amsfonts}
\usepackage{amssymb}
\usepackage{graphicx}
\usepackage{psfrag}

\newtheorem{thm}{Theorem}
\newtheorem{lem}[thm]{Lemma}
\newtheorem{cor}[thm]{Corollary}
\newtheorem{prop}[thm]{Proposition}

\theoremstyle{definition}
\newtheorem{df}[thm]{Definition}

\begin{document}

\title[Simply connected minimal symplectic 4-manifolds]{Simply connected minimal symplectic 4-manifolds\\ with signature less than  --1 }

\author[A. Akhmedov]{Anar Akhmedov}
\author[S. Baldridge]{Scott Baldridge}
\author[R. \.{I}. Baykur]{R. \.{I}nan\c{c} Baykur}
\author[P. Kirk]{Paul Kirk}
\author[B. D. Park]{B. Doug Park}

\date{May 5, 2007}

\address{School of Mathematics,
Georgia Institute of Technology \newline \hspace*{.375in} Atlanta,
GA, 30322, USA} \email{\rm{ahmadov@math.gatech.edu}}

\thanks{A. A. was partially supported by the NSF grant FRG-0244663.  S. B.  was partially supported by the NSF  grant DMS-0507857. R. \.{I}. B. was partially supported by the NSF Grant DMS-0305818. P. K.  was partially supported by the NSF grant DMS-0604310. B.D.P. was partially supported by an NSERC discovery grant.}

\address{Department of Mathematics, Louisiana State University \newline
\hspace*{.375in} Baton Rouge, LA 70817, USA}
\email{\rm{sbaldrid@math.lsu.edu}}

\address{Department of Mathematics, Michigan State University \newline
\hspace*{.375in} East Lansing MI 48824, USA}
\email{\rm{baykur@msu.edu}}

\address{Mathematics Department, Indiana University \newline
\hspace*{.375in} Bloomington, IN 47405, USA}
\email{\rm{pkirk@indiana.edu}}

\address{Department of Pure Mathematics, University of Waterloo\newline
\hspace*{.375in} Waterloo, ON, N2L 3G1, Canada}
\email{\rm{bdpark@math.uwaterloo.ca}}


\maketitle

\begin{abstract} For each pair $(e,\sigma)$  of integers satisfying $2e+3\sigma\ge 0$, $\sigma\leq -2$, and $e+\sigma\equiv 0\pmod{4}$, with four exceptions,  we construct a minimal, simply connected symplectic 4-manifold
with Euler characteristic $e$ and signature $\sigma$.  We also produce simply connected, minimal symplectic 4-manifolds
with signature zero (resp. signature $-1$)   with Euler
characteristic  $4k$ (resp. $4k+1$) for all $k\ge 46$ (resp. $k\ge
49$). \end{abstract}

%**************************************************
%******* Introduction *****************************
%**************************************************

\section{Introduction}

 In \cite{BK4},  a closed, simply connected, minimal symplectic 4-manifold with Euler characteristic $6$ and signature $-2$ is constructed.  This manifold contains a symplectic genus 2 surface with trivial normal bundle and simply connected complement and also contains  two Lagrangian tori with special properties. In this article we use this manifold and apply  standard  constructions to fill out the part of the  symplectic geography plane corresponding to signature less than $-1$.  Recall that Taubes proved (\cite{taubes3, taubes1,taubes2}, also Li-Liu \cite{LiL}) that minimal simply connected symplectic 4-manifolds satisfy $2e+3\sigma\ge 0$, where $e$ denotes the Euler characteristic and $\sigma$ the signature.  Moreover, every symplectic 4-manifold satisfies $e+\sigma\equiv 0\pmod{4}.$

 Our main result is the following.

\noindent{\bf Theorem A.} {\em  Let $\sigma$ and $e$ denote integers satisfying
 $2e+3\sigma\ge 0,$ and $e+\sigma\equiv 0\pmod{4}.$

\noindent If, in addition, $$\sigma\leq -2,$$
 then there exists a simply connected minimal symplectic $4$-manifold with signature $\sigma$ and  Euler characteristic $e$ and odd intersection form, except possibly for $(\sigma,e)$ equal to $(-3, 7)$, $(-3,11)$, $(-5,13)$, or $(-7,15)$.}

\bigskip

In terms of $c_1^2$ and $\chi_h$,   we construct symplectic manifolds realizing all pairs of integers satisfying $0\leq c_1^2\leq 8\chi_h-2$ except $(c_1^2, \chi_h) = (5,1)$, $(13, 2)$, $(11,2)$, and $(9,2)$.

\medskip
Using Freedman's theorem \cite{Freedman}  and Taubes's results \cite{taubes1,taubes2} this theorem can be restated by saying that there exists a minimal symplectic manifold homeomorphic but not diffeomorphic to $m \CP^2\# n\bCP^2$ whenever $m+2\leq n\leq 5m+4$ and $m$ is odd, except possibly for $(m,n)= (1,4),(3, 6), (3,8),$ or $(3,10)$.  The existence of minimal symplectic 4-manifolds homeomorphic to $m\CP^2\# n\bCP^2$ for these four pairs remains  an open problem (as far as we know).

\medskip

The study of the geography problem for symplectic 4-manifolds has been an area of active study in recent years. Jongil Park pioneered the problem of systematically filling in large regions of the geography plane in a series of articles \cite{JP1, JP2, JP6} (the smooth geography
also has a long history, cf. \cite{FS8}). More recent articles have focused on the problem of constructing  small symplectic examples; this includes the articles \cite{JP5, FS7, SS, A, ABP, AP, BK2, BK4, FSP}.

\bigskip
The methods in this article are based on inductive constructions to produce simply connected manifolds starting with a few basic models. Although there are some formal similarities between some of the fundamental group calculations carried out in this article and those in the articles \cite{A, BK2, BK4, AP, ABP,FSP}, there is an important difference between the methods used in those articles and the methods of the present article, as we now explain.

In those articles, the mechanism used to kill fundamental groups comes down to  establishing precise enough control over a group presentation to conclude that all generators die.  This is a  subtle process which depends critically on properly identifying words in fundamental groups, since e.g. in a group a pair of elements $x,y$ might commute but their conjugates $gxg^{-1}, hyh^{-1}$ need not.

By contrast, the mechanism of the present paper is much softer. We use standard symplectic constructions pioneered by Gompf \cite{Gompf}  and Luttinger \cite{Lut} to kill a generator outright; subsequent generators then are killed by a simple argument. In particular,  although we are explicit and careful in our fundamental group calculations in Theorem \ref{thm:technical_lem}, Lemma \ref{gen2xgen2},  and elsewhere, the reader  will quickly understand that our results  follow  as easily if one  only knows the statements  up to conjugacy.

To illustrate this point, in the statement of Lemma \ref{gen2xgen2}, the expressions for $\mu_6, m_6, \ell_6$ are long, but it is straightforward to see    that, {\em up to conjugacy}, $\mu_6=[a_1,x_2], m_6=y_2, \ell_6=b_1^{-1}$. This less precise information is quite sufficient to prove the results of this article.

\bigskip

The construction is also suitable to  fill out a large region of the geography
plane starting with any given symplectic 4-manifold with given
characteristic numbers and containing a square zero symplectic torus. For example, Theorem
\ref{wedge}  roughly says that given a symplectic 4-manifold $X$,
one can construct a new symplectic manifold $Y$ with the same
fundamental group as $X$ and satisfying $c_1^2(Y)=c_1^2(X)+c$ and
$\chi_h(Y)=\chi_h(X) +\chi$, for any $(c,\chi)$ in the cone $0\leq
c\leq 8\chi-2$.

Since it is known how to produce manifolds with positive signature (\cite{stip}) we apply this result to a positive signature symplectic 4-manifold  and prove the following (see Theorem \ref{positive}).

\medskip

\noindent{\bf Theorem B.} {\em
 For all integers $k\ge 45$, there exists a simply connected minimal symplectic $4$-manifold $X_{2k+1, 2k+1}$ with Euler characteristic $e=4k+4$ and signature $\sigma=0$.

 For all integers $k\ge 49$ there exists a simply connected minimal symplectic $4$-manifold $X_{2k-1,2k}$
with Euler characteristic $e=4k+1$ and signature $\sigma=-1$.
}

\medskip
All the manifolds we produce have odd intersection forms.
Hence there remain 4 minimal simply connected symplectic odd 4-manifolds of signature less than or equal to $-2$, 97 minimal simply connected symplectic odd 4-manifolds of non-positive signature, and roughly 280 minimal simply connected odd symplectic 4-manifolds of signature less than or equal to 4  yet to be constructed.

\bigskip

We finish this introduction  with a brief description of the proofs. We   start with  three models, the minimal symplectic 4-manifolds $B,C,D$. These manifolds have Euler characteristic $6,8,$ and $10$ and signatures $-2,-4, $ and $-6$ respectively. Each contains a disjoint pair of homologically independent Lagrangian  tori $T_1$ and $T_2$ with nullhomotopic meridians and whose complement has fundamental group $\ZZ\oplus \ZZ$.  Moreover, $\pm 1$ Luttinger surgery (see Section \ref{lutsurgsec}) along certain curves  on one or both of these tori yields a minimal symplectic 4-manifold.

We then produce a family $B_g, g\in \ZZ$, of minimal symplectic 4-manifolds with  Euler characteristic  $6+4g$ and signature $-2$ by taking a symplectic sum of $B$ with a minimal manifold constructed from Luttinger surgery on a product of surfaces. This family $B_g$ again contains a pair of Lagrangian tori $T_1, T_2$ with the same properties as those in $B,C,D$.

Taking the symplectic sum of many copies of $B,B_g, C,D$  (and, if needed, the elliptic surfaces $E(k)$)  along their tori   and performing a  $+1$ Luttinger surgery on the unused Lagrangian tori yields   our even signature examples. Showing that the fundamental group vanishes is simple since the fundamental groups of $B,B_g, C,D$ and the homomorphisms induced by the inclusions of the tori are known.  Usher's theorem \cite{usher} easily implies that the result is minimal.  The manifolds $B,B_g, C,D$ contain $-1$ surfaces disjoint from  the $T_i$ which survive to $-1$ surfaces in the symplectic sum and  hence the result has an odd intersection form.

Producing odd signature manifolds follows the same general approach, but requires several small model manifolds with appropriate Lagrangian tori to use as seeds for the symplectic sums.
   The construction is not quite as clean as in the even signature case.

   We construct  a minimal symplectic 4-manifold  $P_{5,8}$ with fundamental group $\ZZ$, Euler characteristic 15, and signature $-3$. We also use the minimal simply connected symplectic 4-manifold $S_{1,1}$ constructed by Gompf with Euler characteristic $23$ and signature $-15$.  These three manifolds each contain  a  Lagrangian or symplectic torus appropriate for taking   symplectic sums with many copies of $B, B_g,C,$ and $D$, and as in the even case this produces minimal simply connected 4-manifolds of odd signature less than or equal to $-5$.

The signature $-3$ examples are constructed by a separate argument, and a few small examples  not covered by our general construction are  culled from the literature (i.e. $ (\sigma,e)=(-7,11)$, $(-13,21)$,
$(-11,19)$, $(-5,9)$) or constructed explicitly
($(\sigma,e)=(-5,17)$, $(-7,19)$, $(-9,21)$).

\bigskip

The authors would like to thank R.E. Gompf for helpful comments.

\section{Luttinger surgery}\label{lutsurgsec}  Given any Lagrangian torus $T$ in a
symplectic 4-manifold $M$, the Darboux-Weinstein theorem
 \cite{MS} implies that there is a parameterization of a tubular neighborhood of $T^2\times D^2\to nbd(T)\subset M$ such that the image of $T^2\times \{d\}$ is Lagrangian for all $d\in D^2$. Choosing any point $d\ne 0$ in $D^2$ gives a push off $F_d:T\to T^2\times \{d\}\subset M-T$ called the {\em Lagrangian push off} or {\em Lagrangian framing}.  Given any embedded curve $\gamma\subset T$, its image $F_d(\gamma)$ is called the {\em Lagrangian push off} of $\gamma$.

 Any curve isotopic to  $\{t\}\times \partial D^2\subset \partial (nbd(T))$ will be called a {\em meridian} of $T$ and typically denoted by $\mu_T$.  In this article we will typically fix a pair of embedded curves on $T$ intersecting transversally in one point and denote the two Lagrangian push offs by $m_T$ and $\ell_T$. The triple $\mu_T, m_T ,\ell_T$ generate
   $H_1(\partial(nbd(T)))$. Since the 3-torus has abelian fundamental group we may choose a base point $t$ on $\partial(nbd(T))$ and unambiguously refer to $\mu_T, m_T, \ell_T \in\pi_1(\partial(nbd(T)),t)$.

 The push offs and meridian  are used to specify coordinates for a {\em  Luttinger surgery}. This is   the process of removing a tubular neighborhood of $T$ in $M$ and re-gluing it  so that the embedded curve representing $\mu_T  m_T^p \ell_T^q$  bounds a disk for some pair of integers $p,q$.  The   resulting 4-manifold admits a symplectic structure whose  symplectic form  is unchanged away from a neighborhood of $T$   (\cite{Lut, ADK}).

When the base point $x$ of $M$ is chosen off the boundary of the tubular neighborhood of $T$, the based loops $\mu_T, m_T,$ and $\ell_T$ are to be joined to $x$ {\em by the same path} in $M-T$. These curves then define elements of $\pi_1(M-T,x)$. With $p,q$ as above, the 4-manifold resulting from    torus surgery on $M$ has fundamental group
$$\pi_1(M-T,x)/N(\mu_T  m_T^p \ell_T^q),$$
where $N(\mu_T  m_T^p \ell_T^q)$ denotes the normal subgroup generated by $\mu_T  m_T^p \ell_T^q$.

We will only need the cases $(p,q)=(\pm 1,0)$ or $(0, \pm 1)$ in this article, i.e. {\em $\pm 1$ Luttinger surgery along $m_T$} or $\ell_{T}$.

\section{The fundamental group of the complement of tori in the
product of surfaces}\label{HxK}

Let $F$ be a genus $f$ surface, with $f\ge 2$. Choose a base point $h
$ on $F$ and   pairs $x_i,y_i, i=1,\dots, f$   of circles forming a
symplectic basis, with $x_i,y_i$ intersecting at $h_i\in F$.  Choose
paths $\alpha_i$ from $h$ to $h_i$, so that the loops
$$\tx_i=\al_i x_i\al_i^{-1}  \text{ and } \ty_i=\al_iy_i\al_i^{-1}$$
generate $\pi_1(F, h)$. Let $Y_i$ be a circle parallel to $y_i$ which
misses $\alpha_i$.

Let  $G$ a genus $g$ surface.  Choose a   base point $k$ on $G$, and
$g$ pairs $a_1,b_1,\cdots, a_g,b_g$   of circles forming a symplectic
basis, with $a_i, b_i$ intersecting at $k_i$.
Choose paths $\be_i$ from $h$ to $h_i$, so that the loops
$$\ta_i=\be_ia_i\be_i^{-1} \text{ and } \tb_i=\be_ib_i\be_i^{-1}$$
generate $\pi_1(G,k)$.  Choose parallel copies $A_i$ of $a_i$ and $B_i
$ of $b_i$ which miss the paths $\beta_i$. In Figure 1 we illustrate
the notation when   $f=2$ and $g=3$.

\begin{figure}[h]
\begin{center}
\small
\psfrag{y1}{$y_1$} \psfrag{x1}{$x_1$} \psfrag{x2}{$x_2$} \psfrag{y2}
{$y_2$}
\psfrag{X}{$X$} \psfrag{h}{$h$}\psfrag{h1}{$h'$}\psfrag{K}{$K$}\psfrag
{H}{$H$}
\psfrag{Y1}{$Y_1$} \psfrag{A1}{$A_1$} \psfrag{A2}{$A_2$} \psfrag{Y2}
{$Y_2$}\psfrag{A4}{$A_4$}\psfrag{a1}{$a_1$}\psfrag{a2}{$a_2$}\psfrag
{b1}{$b_1$}
\psfrag{b2}{$b_2$}\psfrag{a3}{$a_3$}\psfrag{b3}{$b_3$}\psfrag{al1}{$
\alpha_1$}\psfrag{al2}{$\alpha_2$}
\psfrag{be1}{$\beta_1$}\psfrag{be2}{$\beta_2$}\psfrag{be3}{$\beta_3$}
\psfrag{B1}{$B_1$} \psfrag{k}{$k$}\psfrag{u}{$u$}\psfrag{k1}{$k_1$}
\psfrag{k2}{$k_2$}
\includegraphics[scale=.80]{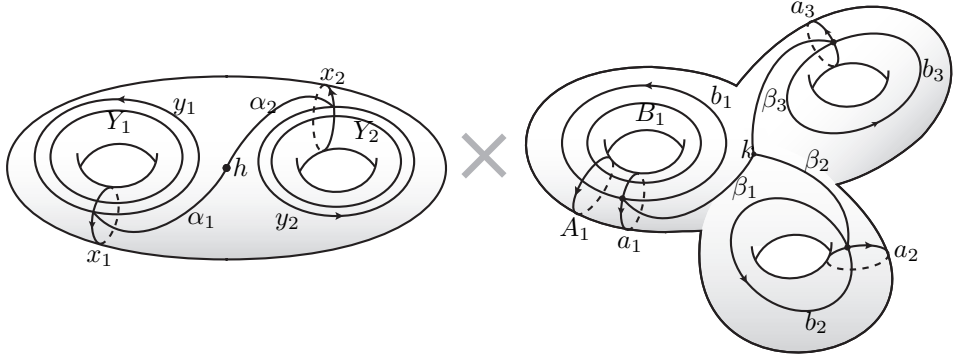}
\caption{The surface $F \times G$.}
\end{center}
\end{figure}

The product $F\times G$  contains the union of the two symplectic
surfaces $F\times\{k\}\cup \{h\}\times G$ meeting at  $(h,k)$.
There is an identification $\pi_1(F\times G, (h,k))=\pi_1(F,h)\times
\pi_1(G,k)$ which associates the loop $\tx_i\times \{k\}$ to $(\tx_i,
1)$, $\ty_i\times\{k\}$ to $(\ty_i,1)$, $\{h\}\times \ta_i$ to $(1,
\ta_i)$ and $\{h\}\times \tb_i$ to $(1,\tb_i)$. In other words, the
homomorphisms induced by the inclusions $F\times\{k\}\subset F\times G
$ and $\{h\}\times G\subset F\times G$ present $\pi_1(F\times G,
(h,k))$ as the product of $\pi_1(F,h)$ and $ \pi_1(G,k)$.

When there is no chance of confusion we denote the $2f+2g$ loops $
\tx_i\times\{k\},\ty_i\times\{k\},\{h\}\times \ta_i,\{h\}\times \tb_i
$  simply by $\tx_i,\ty_i,\ta_i,\tb_i$.  These are loops in $F\times G
$ based at $(h,k)$.

\medskip

The product $F\times G$ contains $2g$   Lagrangian tori
$$Y_1\times A_j, \ \  Y_2\times B_j, \ \ \ j=1,\dots, g.$$
These $2g$ tori are pairwise disjoint and miss $(F\times\{k\})\cup (\{h\}
\times G)$.

Let $N$ denote a tubular neighborhood of the union of these  $2g$ tori:
$$N=nbd\big((\cup_{j} Y_1\times A_j) \cup (\cup_{j}Y_2\times B_j)\big)\subset F\times
G.$$
The loops $\tx_i,\ty_i,\ta_i,\tb_i$ are loops
in $F\times G-N$ based at $(h,k)$.

Typically, removing a surface
from a 4-manifold increases the number of generators of  the
fundamental group, but since these tori respect the product structure
one can prove the following theorem.

\begin{thm} \label{thm:technical_lem} The $2f+2g$ loops
$\tx_1,\ty_1,\dots ,\tx_f,\ty_f, \ta_1,\tb_1,\cdots, \ta_g,\tb_g$
generate
$\pi_1(F\times G -N, (h,k))$.

Moreover, there are paths $d_{j}:[0,1]\to F\times G -N$  from $(h,k)$ to
the boundary of the tubular neighborhood of $Y_1\times A_j$ and $e_
{j}:[0,1]\to F\times G -N$  from $(h,k)$ to the boundary of the tubular
neighborhood of $Y_2\times B_j$ so that with respect to these paths,
the meridian and two Lagrangian push offs of $Y_1\times A_j$ are
homotopic in $F\times G -N$ rel endpoint
to
$$\mu_{Y_1\times A_j}=[\tx_1,\tb_j], \ \ m_{Y_1\times A_j}= \ty_1,\ \  \ell_
{Y_1\times A_j}=\ta_j$$
and
the meridian and two Lagrangian push offs of $Y_2\times B_j$ are
homotopic in $F\times G -N$ rel endpoint
to
$$\mu_{Y_2\times B_j}=[\tx_2,\ta_j],\ \  m_{Y_2\times B_j}= \ty_2, \ \ \ell_
{Y_2\times B_j}=\tb_j.$$
\end{thm}
\begin{proof}

Before we start the proof, we give an indication of how it will
proceed. Note that   $\cup_j(Y_1\times A_j)=Y_1\times (\cup_j A_j)$
lies on $Y_1\times G$ and that $Y_2\times (\cup_j B_j)$ lies on $Y_2
\times G$. Thus  $F\times G -N$ can be constructed by cutting $F
\times G$ along the hypersurface $(Y_1\cup Y_2)\times G$, and then
regluing the two copies of $Y_1\times G$ only along the complement of
a neighborhood of the $A_i$ and regluing the two copies of $Y_2\times
G$ only along the complement of a neighborhood of the $B_i$.
However, in order to use the Seifert-Van Kampen theorem, the subsets
and their intersection in a decomposition are required to be
connected, and so we need to modify the  decomposition slightly.

\bigskip

Let $P_1$ be the annulus in $F$ bounded by $y_1$ and $Y_1$. Similarly
let $P_2$ denote the annulus  in $F$ bounded by $y_2$ and $Y_2$.  Let
$\al$ denote the arc $(\al_1\cup \al_2)\times\{k\}$.  Let $\ga_1$
denote the arc $(x_1\cap P_1)\times \{k\}$; it spans the two circles
$y_1$ and $Y_1$. Similarly let $\ga_2$ denote the arc $(x_2\cap P_2)
\times \{k\}$.  See Figure 2.

\begin{figure}[h]
\psfrag{al}{$\alpha$}\psfrag{P1}{$P_1$}\psfrag{y2}{$\gamma_2$}\psfrag
{P2}{$P_2$}\psfrag{y1}{$\gamma_1$}
\begin{center}
\small
\includegraphics[scale=1.00]{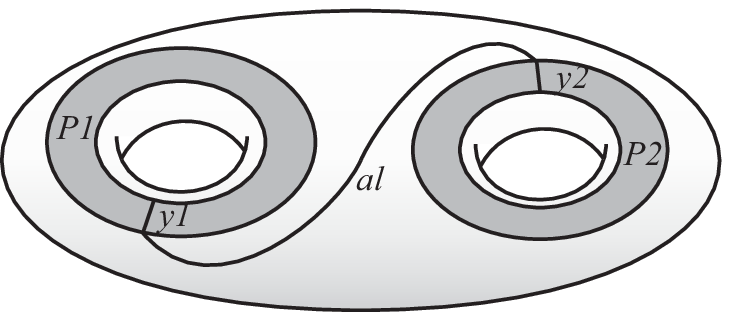}
\caption{   }
\end{center}
\end{figure}

\noindent Set
$$S_1=(P_1\times G)\cup \al\cup (P_2\times G)$$
and
$$S_2= (\ga_1\cup\ga_2)\cup\big( (F-int(P_1\cup P_2))\times  G\big) .$$
Then in $F\times  G$,
$S_1\cap S_2$ is the union of four copies of $S^1\times  G$ together
with three arcs which connect the four components. In particular,
$S_1, S_2$ and $S_1\cap S_2$ are connected and contain the base point
$(h,k)$.

Let $ G_A= G-nbd(\cup_iA_i)$ denote the complement of an open tubular
neighborhood of the $A_i$ in $ G$. Since  the $A_i$ do not disconnect
$ G$, $ G_A$ is path connected.
Similarly let $ G_B= G - nbd(\cup_i B_i)$  denote the complement of
an open tubular neighborhood of the $B_i$ in $ G$.

To construct $F\times  G -N$ we form the identification space
$$ F\times  G -N= S_1\sqcup S_2/\sim$$
by identifying $(f,s)\in S_1$ with its corresponding point $(f',s')$
in $S_2$ {\em  except if }
$f\in Y_1$ and  $s\in nbd(A_i)$ or $f\in Y_2$ and $s\in nbd(B_i)$. In
other words, along $Y_1\times  G$ we identify only the two copies of
$Y_1\times  G_A$ and along $Y_2\times  G$ we identify only the two
copies of $Y_2\times  G_B$.

Hence we have exhibited $F\times  G -N$ as the union of $S_1$ and
$S_2$ with connected intersection
$$S_1\cap S_2=(Y_1\times G_A)\cup \ga_1\cup( y_1\times  G)\cup \al
\cup ( y_2\times  G)\cup\ga_2
\cup(Y_2\times  G_B).$$

It is easy to see that $\pi_1(S_1\cap S_2, (h,k))\to \pi_1(S_1, (h,k))
$ is surjective.  Indeed, one can use the product parameter in the
annuli $P_1$ and $P_2$ to define a deformation retract (fixing $\al$
and hence also $(h,k)$) of $S_1$ to the subset $( y_1\times  G)\cup
\al\cup ( y_2\times  G)$ of $S_1\cap S_2$.

The Seifert-Van Kampen theorem applies and implies that there is a
surjection
$$\pi_1(S_2,(h,k))\to \pi_1(F\times  G -N,(h,k))$$
induced by   inclusion.

We will show that the image of
$\pi_1(S_2,(h,k))\to  \pi_1(F\times  G -N,(h,k))$  is generated by
the  loops
   $\tx_1, \ty_1,\tx_2,\ty_2,\cdots, \tx_f,\ty_f, \ta_1,\tb_1,\cdots,
\ta_g, \tb_g$. Notice that all these loops are contained in $S_2$.

   We find generators for  $\pi_1(S_2,(h,k))$. This is again a
straightforward application of the Seifert-Van Kampen theorem, as we
will now show.

Since the arcs $\ga_i$ are just segments that lie on $\tx_i$ (and the
rest of the loops $\tx_i$ lie in $S_2$), we can decompose $S_2$ as
$$S_2=(\tx_1\cup \tx_2)\cup\big( (F-int(P_1\cup P_2))\times G\big) .$$
The intersection of the two pieces in this decomposition is the
(contractible) set $\tx_1\cup\tx_2-(\ga_1\cup \ga_2)$.

Hence $\pi_1(S_2,(h,k))$ is generated by $\tx_1,\tx_2$ and any set of
generators of
$$\pi_1((F-int(P_1\cup P_2))\times G, (h,k))=\pi_1(F-int(P_1\cup
P_2)),h)\times\pi_1( G, k).$$
The loops $\ta_i,\tb_i$ generate $\pi_1( G)$.  The space $ F-int(P_1
\cup P_2))$ is a 4-punctured genus $f-2$ surface. Its fundamental
group is generated by $\ty_1,\ty_2, \tx_3,\ty_3,\cdots,\tx_f,\ty_f$
and one other loop $\tau$ based at $h$ which is obtained by traveling
from the base point to a point on the boundary component $Y_1$,
following $Y_1$, then returning to the base point.

We have shown that the loops $\tx_1, \ty_1,$ $\tx_2,\ty_2,$ $\tx_3,\ty_3,$ $
\cdots,\tx_f,\ty_f,$ $\ta_1,\tb_1,\cdots, \ta_g, \tb_g$, and $\tau\times
\{k\}$ generate  $\pi_1(S_2,(h,k))$. Hence, considered as loops in
$F\times  G -N$,  they generate $   \pi_1(F\times  G -N,(h,k))$.  We
need only show that the generator $\tau\times\{k\}$ is not needed.
But this is obvious since
   $\tx_1,\ty_1,\cdots,\tx_f,\ty_f $ and $\tau\times\{k\}$ all lie on the
surface
   $F\times \{k\}\subset F\times  G -N$,  and $\tx_1,\ty_1,\cdots, \tx_f,
\ty_f $ generate $\pi_1(F\times \{k\}, k)$.

\medskip

We next turn to the problem of expressing the meridians and
Lagrangian push offs of the generators of the Lagrangian tori $Y_1
\times A_j$, $Y_2\times B_j$ in terms of the loops  $\tx_1, \ty_1,
\cdots, \tx_f,\ty_f,$ $\ta_1,\tb_1,\cdots, \ta_g, \tb_g$. We do this for
$Y_1\times A_1$.  Symmetric arguments provide the analogous
calculations for the rest.

In Figure 1,  denote by $h_1$ the intersection of $x_1$ and $y_1$
(i.e. the endpoint of $\al_1$) and denote by $k_1$ the intersection
of $a_1$ and $b_1$.  Then the point $(h_1,k_1)$ lies on the boundary
of the tubular neighborhood of $Y_1\times A_1$.

Since we take the product symplectic form on $F\times  G$, referring
to Figure 1 one sees that the loops $y_1\times\{k_1\}$ and $\{h_1\}
\times a_1$ are Lagrangian  push offs of  two generators of $\pi_1(Y_1
\times A_1)$ to the boundary of the tubular neighborhood of $Y_1
\times A_1$.

There is a map of a square into
$F\times  G -N$ given by $\al_1\times\be_1$:
$$\al_1\times\beta_1:[0,1]\times[0,1]\to F\times  G -N.$$
The point $(0,0)$ is mapped to the base point $(h,k)$ of $F\times  G -
N$, and the point $(1,1)$ is mapped to $(h_1,k_1)$.  Thus the
diagonal path $d(t)= (\alpha_1(t),\beta_1(t))$  connects the base
point to the boundary of the tubular neighborhood of $T_1$.

Conjugating by $d$ expresses the Lagrangian push offs as based curves in
$F\times  G -N$.  So
$$m_{Y_1\times A_1}= d*(y_1\times \{k_1\})*d^{-1}\text{ and } \ell_
{Y_1\times A_1}= d*( \{h_1\}\times a_1)*d^{-1}.$$
But $m_{Y_1\times A_1}$ is homotopic rel basepoint in $F\times  G -N$
to $\ty_1$. An explicit homotopy is given by the formula:
$$(s,t)\mapsto
\begin{cases} \big(\al_1(3t),  \be_1((1-s)3t)\big)& \text{ if } 0\leq
t\leq \frac{1}{3}\\
         \big( y_1(3t-1),  \beta_1(1-s)\big)& \text{ if } \frac{1}{3}
\leq t\leq \frac{2}{3}\\
         \big(\al_1(3-3t), \be_1((1-s)(3-3t))\big)& \text{ if } \frac
{3}{3}\leq t\leq 1.
\end{cases}
$$
A similar homotopy, but exchanging the roles of $\al_1$ and $\be_1$,
establishes that
   $\ell_{Y_1\times A_1}$ is homotopic rel basepoint in $F\times  G -N
$ to $\ta_1$.
(These homotopies clearly miss all the other $Y_1\times A_j$ and $Y_2
\times B_j$.)

It remains to calculate the meridian of $Y_1\times A_1$.  For this,
consider the map
$x_1\times b_1:[0,1]\times[0,1]\to F\times  G$. This has image a
torus  intersecting  $T_1$ transversally in one point (near the point
$(x_1(.9),b_1(.9))$, as one sees from Figure 1).  Since
$$(h_1,k_1)=(x_1\times b_1)(0,0)=(x_1\times b_1)(0,1)=(x_1\times b_1)
(1,0)=(x_1\times b_1)(1,1),$$
by conjugating the path that follows the boundary of this square by
the path $d$ from  the  base point $(h,k)$ to $(h_1,k_1)$,  we see
that the meridian $\mu_{Y_1\times A_1}$ is given by the composite
$$\mu_{Y_1\times A_1}= d*(x_1\times \{k_1\})* (\{h_1\}\times b_1)*(x_1
\times \{k_1\})^{-1}*(\{h_1\}\times b_1)^{-1}*d^{-1}.$$

Now $d*(x_1\times\{k_1\})*d^{-1}$ is homotopic rel basepoint to $\tx_1
$ in $F\times  G -N$ by the same argument given above. The key
observation is that   $\beta_1$ misses $A_i$ and $B_i$ for all $i$.
Similarly
$d*(\{h_1\}\times b_1)*d^{-1}$  is homotopic rel basepoint to $\tb_1$
in $F\times  G -N$.  Thus
$$\mu_{Y_1\times A_1}\sim \tx_1* \tb_1*\tx_1^{-1}*\tb_1^{-1}=[\tx_1,
\tb_1].$$

Similar calculations establish all other assertions.
\end{proof}

\section{Telescoping triples and symplectic sums}

Our construction of symplectic 4-manifolds which fill large regions in the geography plane is based on using telescoping symplectic sums along symplectic tori as well as Luttinger surgeries. The basic models in our construction have a convenient property preserved under appropriate symplectic sum, and so  we formalize the property in the following definition.

\begin{df} An ordered  triple  $(X,T_1, T_2)$ where $X$ is a symplectic 4-manifold and
$T_1,T_2$ are disjointly embedded Lagrangian tori is called a {\em telescoping triple} if
\begin{enumerate}
\item The tori $T_1,T_2$ span a 2-dimensional subspace of $H_2(X;\RR)$.
\item $\pi_1(X)\cong \ZZ^2$ and the inclusion induces an isomorphism $\pi_1(X-(T_1\cup T_2))\to \pi_1(X)$ (in particular the meridians of the $T_i$ are trivial in $\pi_1(X-(T_1\cup T_2))$).
\item The image of the homomorphism induced by inclusion $\pi_1(T_1)\to \pi_1(X)$ is a {\em summand} $\ZZ\subset \pi_1(X)$.
\item The homomorphism induced by inclusion $\pi_1(T_2)\to \pi_1(X)$ is an isomorphism.
 \end{enumerate}
 If $X$ is minimal we call $(X,T_1,T_2)$ a {\em minimal telescoping triple}.
 \end{df}

 Note that the order of $(T_1,T_2)$,    matters  in this definition. Notice  also that since the meridians $\mu_{T_1},\mu_{T_2}\in \pi_1(X-(T_1\cup T_2))$ are trivial and the relevant fundamental groups are abelian, the push off of an oriented loop $ \gamma\subset T_i$ into $X-(T_1\cup T_2)$ with respect to any framing of the normal bundle of $T_i$ (e.g. the Lagrangian framing) represents a well defined element of $\pi_1(X-(T_1\cup T_2))$, independent of the choice of framing (and basing).

 The definition of a telescoping triple includes the hypothesis that the Lagrangian tori $T_1$ and $T_2$ are linearly independent in $H_2(X;\RR)$. This implies (\cite{Gompf}) that the symplectic form on $X$ can be slightly perturbed so that one of the $T_i$ remains Lagrangian and the other becomes symplectic. It can also be perturbed so that both become symplectic. Moreover, if $F$ is a symplectic surface in $X$ disjoint from $T_1$ and $T_2$, the perturbed symplectic form can be chosen so that $F$ remains symplectic.
 \bigskip

 Recall that the symplectic sum (\cite{Gompf}) of two symplectic 4-manifolds $X$ and $X'$ along genus $g$ symplectic surfaces $F\subset X$ and $F'\subset X$ of opposite square is a symplectic 4-manifold described topologically as the union
 $$X\#_{F,F'}X'=(X-nbd(F))\cup (X'-nbd(F'))$$
 where the boundaries of the tubular neighborhood are identified by a fiber-preserving diffeomorphism of the corresponding circle bundles. When the surfaces are clear from context we write $X\#_s X'$.

 \begin{prop}\label{symp_sum_type_I_or_II}  Let $(X, T_1,T_2)$  and  $(X', T_1',T_2')$ be two telescoping triples.  Then  for an appropriate gluing map  the triple
 $$(X\#_{T_2,T_1'}X',T_1,T_2') $$
 is again a telescoping triple.

The Euler characteristic and signature of $X\#_{T_2,T_1'}X'$ are given by $e(X)+e(X')$ and $\sigma(X)+\sigma(X')$.
\end{prop}

\begin{proof}     Let $i_j:\pi_1(T_j)\to \pi_1(X)$ be the homomorphisms induced by inclusion  for $j=1,2$.
Choose $x_1, y_1\in \pi_1(T_1)$ so that $x_1$ spans the kernel of  $i_1$ and $i_1(y_1)$ spans the image of $i_1$.  Denote $i_1(y_1)$ by $t$ and choose $s\in \pi_1(X)$ so that $s,t$ forms a basis of $\pi_1(X)$. Then choose generators $x_2,y_2$ for $\pi_1(T_2)$ so that $i_2(x_2)=s$ and $i_2(y_2)=t$.
Thus the inclusions induce
$$x_1\mapsto 1, y_1\mapsto t, x_2\mapsto s, y_2\mapsto t.$$

 Similarly, construct generators $x_1',y_1'$ for $\pi_1(T_1')$, $x_2', y_2'$ for $\pi_1(T_2')$ and $s',t'$ for $\pi_1(X')$.

The inclusion  induces an  isomorphism  $\pi_1(X-(T_1\cup T_2))\to
\pi_1(X)$ and  the boundary of the tubular neighborhood of $T_1$
(resp. $T_2$) is a 3-torus whose fundamental group is spanned by
$x_1,y_1,\mu_{T_1}$ (resp. $x_2,y_2,\mu_{T_2}$) (for definitiveness use the
Lagrangian framing to push the $x_i, y_i$ into the boundary of the tubular neighborhood). Similar assertions hold for $(X',T_1', T_2')$.
The symplectic sum of $X$ and $X'$ along the surfaces $T_2\subset X,
T_1'\subset X'$ can be formed so that the ordered triple
$(x_2,y_2,\mu_2)$ is sent to $(x_1',y_1',\mu_1')$ by the gluing
diffeomorphism (perhaps after a change of orientation on some of the
loops to ensure that the gluing diffeomorphism is orientation
preserving).

The Seifert-Van Kampen theorem and the fact that all meridians are trivial imply that
$$\pi_1(X\#_{T_2,T_1'}X') =\langle s,t,s',t'\ | \ [s,t],[s',t'], s, t(t')^{-1}\rangle=\ZZ s'\oplus \ZZ t'.$$
The inclusion $T_1\subset X\#_{T_2,T_1'}X'$ induces $x_1\mapsto 1, y_1\mapsto t'$. The inclusion $T_2'\subset X\#_{T_2,T_1'}X'$ induces $x_2'\mapsto s', y_2'\mapsto t'$.   Hence $(X\#_{T_2,T_1'}X',T_1,T_2') $
 is indeed a telescoping triple.

 The assertions about the Euler characteristic and signature are clear.
\end{proof}

\medskip
Since the meridians of the Lagrangian tori are trivial in a telescoping triple, one immediately concludes the following.
 \begin{prop}\label{surger} Let  $(X,T_1, T_2)$ be a telescoping triple. Let $\ell_{T_1}$ be a Lagrangian push off of a curve on $T_1$  and $m_{T_2} $ the Lagrangian push off of a curve on $T_2$ so that $\ell_{T_1}$ and $m_{T_2}$ generate $\pi_1(X)$.

 Then
  the symplectic 4-manifold obtained by performing $+1$ Luttinger surgery on $T_1$ along $\ell_{T_1}$ and $+1$ surgery on $T_2$ along $m_{T_2}$ is simply connected. \qed
 \end{prop}
\medskip

We will have frequent use of the following two results. The first is a criterion given by Usher \cite{usher} to determine when a symplectic 4-manifold is minimal. The second is a useful result of T. J. Li which we will use to verify that the hypotheses in Usher's theorem hold in certain contexts.

\begin{thm}[Usher]\label{Ush} Let  $Z=X_1\#_{F_1=F_2} X_2$ denote the symplectic sum of $X_1$ and $X_2$ along symplectic surfaces  $F_i$ of  positive genus $g$.   Then:
\begin{itemize}
\item[(i)] If either $X_1- F_1$ or $X_2- F_2$
contains an embedded symplectic sphere of square $-1$, then $Z$ is
not minimal.
\item[(ii)] If one of the summands $X_i$ (for definiteness, say
$X_1$) admits the structure of an $S^2$-bundle over a surface of
genus $g$ such that $F_1$ is a section of this fiber bundle, then
$Z$ is minimal if and only if $X_2$ is minimal.
\item[(iii)] In all other cases, $Z$ is minimal.
\end{itemize}
\end{thm}

\bigskip

 Corollary 3 of T.J. Li's article \cite{TJli}  provides a useful method to eliminate the first two cases of Usher's theorem in some contexts.

 \begin{thm}[Li] \label{Li} Let $M$ be a symplectic 4-manifold which is not rational or ruled. Then every smoothly embedded $-1$ sphere is homologous to a symplectic $-1$ curve up to sign. If $M$ is the blow up of a minimal symplectic 4-manifold with $E_1,\cdots, E_n$ represented by exceptional curves, then the $E_i$ are the only classes represented by a smoothly embedded $-1$ sphere, hence any orientation preserving diffeomorphism maps $E_i$ to some $\pm E_j$.
 \end{thm}

\section{The model even signature manifolds}

We will setup an inductive argument by constructing telescoping symplectic sums starting with several basic telescoping triples. Proposition \ref{surger} then applies to produce simply connected 4-manifolds.

\bigskip

To begin with, in  \cite[Theorem 20]{BK4}, a minimal  telescoping triple $(B, T_1, T_2)$  is constructed ($B$ is denoted $B_1$ in that article) so that $B$ contains  a genus 2 surface $F$ with trivial normal bundle, and a geometrically dual symplectic $-1$ torus $H_1$. The tori $T_1, T_2$ miss $F\cup H_1$.  Moreover,  $(B-F,T_1,T_2)$ is also a  telescoping triple. These facts follow immediately from the following theorem, which summarizes the assertions established in \cite{BK4}.

 \begin{thm}\label{manB} There exists a minimal symplectic 4-manifold $B$ containing a pair of homologically essential Lagrangian tori $T_1$ and $T_2$ and a square zero symplectic genus 2 surface $F$ so that $T_1,T_2$, and $F$ are pairwise disjoint,    $\eu(B) = 6 $ and $\sigma(B)=-2$, and
\begin{enumerate}
\item $\pi_1(B-(F\cup T_1\cup T_2))= \ZZ^2$, generated by $t_1$ and $t_2$.
\item The inclusion  $B-(F\cup T_1\cup T_2)\subset B$ induces an isomorphism on fundamental groups. In particular the meridians $\mu_F,\mu_{T_1},\mu_{T_2}$ all vanish in $\pi_1(B-(F\cup T_1\cup T_2))$.
\item The Lagrangian push offs  $m_{T_1},\ell_{T_1}$ of $\pi_1(T_1)$ are sent to $1$ and $t_2$ respectively in the   fundamental group  of $ B-(F\cup T_1\cup T_2)$.
\item The Lagrangian push offs  $m_{T_2},\ell_{T_2}$ of $\pi_1(T_2)$ are sent to $t_1$ and $t_2$ respectively in the   fundamental group  of $ B-(F\cup T_1\cup T_2)$.
\item The push off
 $F\subset   B-(F\cup T_1\cup T_2)$ takes the first three generators of a standard symplectic generating set $\{a_1,b_1,a_2,b_2\}$ for $\pi_1(F)$ to $1$ and the last element to $t_2$.
 \item There exists  a symplectic torus $H_1\subset B$ which intersects $F$ transversally once, which has square $-1$, and the homomorphism $\pi_1(H_1)\to \pi_1(B)$ takes the first generator to $1$ and the second to  $ t_1$. Moreover $H_1$ is disjoint from $T_1$ and $T_2$ (see \cite[Proposition 12, Theorem 20]{BK4}).
\end{enumerate}
  \qed
 \end{thm}

  The following is a restatement of \cite[Theorem 13]{BK4}.  We state it formally since we will have frequent need of it.

  \begin{cor} \label{1-3} The symplectic 4-manifold $X_{1,3}$ obtained from $B$ by $+1$ Luttinger surgery on $T_1$ along $\ell_{T_1}$ and $+1$ Luttinger surgery on $T_2$ along $m_{T_2}$ is a
minimal  symplectic 4-manifold    homeomorphic to $\CP^2\# 3\bCP^2$. It contains a genus 2 symplectic surface of square zero with simply connected complement and a symplectic torus $H_1$ of square $-1$ intersecting $F$ transversally and positively in one point.\qed
  \end{cor}

 \medskip

 \begin{cor}\label{cor6} For each $g\ge 0$ there exists a minimal telescoping triple $(B_g,T_1,T_2)$ satisfying $e(B_g)=6+ 4g$ and $\sigma(B_g)=-2$ and containing a square $-1$ genus $ g+1$ surface disjoint from $T_1\cup T_2$.
 \end{cor}
 \begin{proof} To avoid confusing notation, during this proof we denote the symplectic genus 2 surface in $B$  of Theorem
 \ref{manB} by $F_B$.

Take the product $F\times G$ of a genus $2$ surface $F$ and  a genus $g$ surface
 $G$, as in Section \ref{HxK}.
 Let $Z_g$ denote the 4-manifold obtained from   $F\times G$   by performing  $-1/1$ Luttinger surgeries on the $2g$ disjoint Lagrangian tori $Y_1\times A_i $ and $Y_2\times B_i$ along the curves  $\ell_{Y_1\times A_i }=a_i$ and $\ell_{Y_2\times B_i }=b_i$.  Then by Theorem  \ref{thm:technical_lem}  the fundamental group of $Z_g$ is generated by  loops $\tx_1,\ty_1,\tx_2,\ty_2, \ta_1,\tb_1,\cdots, \ta_g,\tb_g$ and the relations
 $$[\tx_1,\tb_i]=\ta_i, [\tx_2,\ta_i]=\tb_i$$ hold in $\pi_1(Z_g)$. Moreover, the standard symplectic generators for $\pi_1(F)$ are sent to $\tx_1,\ty_1,\tx_2,\ty_2$ in $\pi_1(Z_g)$.

Since the meridian $\mu_{F_B}$ of $F_B\subset B$ is trivial,
the symplectic sum   of $B$ with $Z_g$ along their genus 2 symplectic surfaces $F_B\subset B$ and $F=F\times\{k\}\subset F\times G$
$$B_g=B\#_{F_B, F }Z_g$$
 has fundamental group a quotient of $(\ZZ t_1\oplus \ZZ t_2)*\pi_1(Z_g)$. We choose this symplectic sum so that the generators $a_1,b_1,a_2,b_2$ for $\pi_1(F_B)$ are  identified (in order) with the generators $\tx_1,\ty_1,\tx_2,\ty_2$.

  The fifth assertion of Theorem \ref{manB} shows that $\tx_1,\ty_1,$ and $ \tx_2$ are trivial in
  $\pi_1(B_g)$. The relations coming from the Luttinger surgeries then show  that $\ta_i=1=\tb_i$. Since $b_2=\ty_2$ is identified with $t_2$, $\pi_1(B_g)$ is generated by $t_1 $ and $t_2$. A calculation using the Mayer-Vietoris sequence shows that $H_1(B_g)=\ZZ^2$, and so $\pi_1(B_g)=\ZZ t_1\oplus \ZZ t_2$. Hence $(B_g,T_1,T_2)$ is a telescoping triple, as desired.

The Euler characteristic of $B_g$ is calculated as
$$e(B_g)=e(B)+ e(F\times G)+4=6 + 4g-4 +4= 6+4g,$$ and the signature is computed by Novikov additivity: $\sigma(B_g)=\sigma(B)=-2$.

 The torus $H_1$ in $B-N$ geometrically dual to $F_B$ can be lined up with one of the parallel copies $\{z\}\times G$ in $F\times G-N$ (i.e. take a {\em relative symplectic sum}, \cite{Gompf}) to produce a square $-1$ genus $g+1$ surface in $B_g$.

 Minimality follows from \cite[Lemma 2]{BK4}, which shows that $Z_g$ is minimal (its universal cover is contractible, so $\pi_2(Z_g)=0$) and Usher's theorem, Theorem \ref{Ush}.
 \end{proof}

 \medskip
 We can also produce telescoping triples  with odd signature starting with $B$.   Recall that a symplectic 4-manifold $X$ containing a symplectic surface $F$ is called {\em relatively minimal} if every $-1$ sphere in $X$ intersects $F$.

 \begin{lem}\label{manA}  The blow up $A=B\# \bCP^2$ contains a genus 3 symplectic surface $F_3$ with trivial normal bundle and two Lagrangian tori $T_1$ and $T_2$ so that the surfaces $F_3, T_1,T_2$ are pairwise disjoint,   $(A,F_3)$ is relatively minimal, and:
 \begin{enumerate}
\item $\pi_1(A-(F_3\cup T_1\cup T_2))= \ZZ^2$, generated by $t_1$ and $t_2$.
\item The inclusion  $A-(F\cup T_1\cup T_2)\subset A$ induces an isomorphism on fundamental groups. In particular the meridians $\mu_{F_3},\mu_{T_1},\mu_{T_2}$ all vanish in $\pi_1(A-(F_3\cup T_1\cup T_2))$.
\item The Lagrangian push offs  $m_{T_1},\ell_{T_1}$ of $\pi_1(T_1)$ are sent to $1$ and $t_2$ respectively in the   fundamental group  of $ A-(F_3\cup T_1\cup T_2)$.
\item The Lagrangian push offs  $m_{T_2},\ell_{T_2}$ of $\pi_1(T_2)$ are sent to $t_1$ and $t_2$ respectively in the   fundamental group  of $ A-(F_3\cup T_1\cup T_2)$.
\item There is a standard symplectic generating set  $\{a_1,b_1,a_2,b_2, a_3,b_3\}$for $\pi_1(F_3)$ so that the push off
 $F_3\subset   A-(F_3\cup T_1\cup T_1)$ takes $b_2$ to $t_2$, $b_3$ to $t_1$, and all other generators to $1$.
\end{enumerate}
In particular, $(A,T_1,T_2)$ is a telescoping triple.
 \end{lem}
 \begin{proof} The 4-manifold $B$ of Theorem \ref{manB}  contains a symplectic genus 2 surface $F$ of square zero and a geometrically dual symplectic torus $H_1$ of square $-1$. Symplectically resolve the union $F\cup H_1$ to get $F'_3$, a genus three symplectic surface in $B$ which misses $T_1$ and $T_2$. The surface $F'_3$ has square $(F+H_2)^2= 1$. Blow up $B$ at one point on $F'_3$ to construct $A$ and denote the proper transform of $F_3'$ by $F_3$.

 Since $F_3$ has a geometrically dual 2-sphere (the exceptional sphere), the meridian of $F_3$ in  $A-F_3\subset F_3$ is nullhomotopic.  The rest of the fundamental group assertions follow from Theorem \ref{manB}.

 Although $A$ is not minimal, T.J. Li's theorem  (Theorem \ref{Li})   implies that every $-1$ sphere in $A$ intersects $F_3$, since $B$ is minimal, and neither rational nor ruled.

  \end{proof}

  Note that Luttinger surgery on $T_1$ and $T_2$ in $A$ produces a symplectic 4-manifold homeomorphic to $\CP^2\# 4 \bCP^2$, but this manifold is not minimal; it is just the blow up $X_{1,3}\# \bCP^2$.   We do not know how to produce a minimal symplectic 4-manifold with this homeomorphism type.
  \bigskip

We next produce a 4-manifold $C$   with $e=8$ and $\sigma=-2$ by stopping  the construction of a  minimal symplectic 4-manifold homeomorphic to $\CP^2\# 5\bCP^2$ in the proof of  \cite[Theorem 10]{BK4} before the last 2 Luttinger surgeries to obtain the following.

\begin{thm}\label{manC} There exists a  minimal telescoping triple $(C, T_1,T_2)$ with $\eu(C) = 8$ and $\sign(C)=-4$. Moreover, $C$ contains a square $-1$ torus disjoint from $T_1\cup T_2$.
 \end{thm}

\begin{proof}
We follow the notation and proof of \cite[Theorem 10]{BK4}.  By {\em not} performing the Luttinger surgeries on the tori $T_3$ and $T_4$,    one obtains a minimal symplectic 4-manifold $C$  such that $\pi_1(C - (T_3\cup T_4))$  is generated by the two commuting elements $ y$ and $ a_2$.    The Mayer-Vietoris sequence shows that $H_1(C-(T_3 \cup T_4);\BZ)=\BZ^2$, and so $\pi_1(C-(T_1\cup T_2))=\BZ y \oplus \BZ a_2$.  The  meridians and Lagrangian push offs of $T_3$ and $T_4$ are given by
$\mu_{T_3}=1, m_{T_3}=1, \ell_{T_3}=a_2$ and $\mu_{T_4}=1, m_{T_4}=y, \ell_{T_4}=a_2$.
Thus $(C,T_3, T_4)$ is a telescoping triple.

We relabel $T_3$ by $T_1$ and $T_4$ by $T_2$.

The $-1$ torus comes about from the construction. Briefly, $C$ is obtained by performing Luttinger surgeries on the symplectic sum $(T^2\times F_2)\#_s (T^2\times S^2\#4\bCP^2)$ along the genus 2 surface $\{x\}\times F_2$ in $T^2\times F_2$ and the genus 2 surface $F_2'\subset  (T^2\times S^2)\#4\bCP^2$ obtained by resolving the singularities of $(T^2\times\{p_1\})\cup (\{q\}\times S^2)\cup
(T^2\times\{p_2\})$ and blowing up 4 times at points on this genus 2 surface.  One can choose a square zero torus of the form $T^2\times \{y\}\subset T^2\times F_2$  which matches up (i.e. take a relative symplectic sum) with one of the four exceptional curves to provide  a $-1$ symplectic torus  disjoint from the Lagrangian tori where the Luttinger surgeries are performed.
\end{proof}

The symplectic 4-manifold $X_{1,5}$ obtained from $C$ by $+1$ Luttinger surgeries on $T_1$ and $T_2$ as in Proposition \ref{surger} is   minimal and homeomorphic to $\CP^2\#5\bCP^2$ (\cite{BK4}).

\medskip

Our next small model is a minimal telescoping triple built in the process of  constructing  a minimal symplectic 4-manifold homeomorphic to $\CP^2\# 7\bCP^2$ in \cite[Theorem 8]{BK4}. One stops the construction before performing the  2 Luttinger surgeries, and these unused tori provide the desired $T_1$ and $T_2$.

\begin{thm}\label{BTypeII} There exists a minimal telescoping triple $(D,T_1,T_2)$ with $\eu(D)=10$ and $\sigma(D)=-6$. Moreover, $D$ contains a square $-1$ torus disjoint from $T_1\cup T_2$. \qed
\end{thm}
 \begin{proof} The proof is similar to that of Theorem \ref{manC}. We follow the notation and proof of \cite[Theorem 8]{BK4}. The 4-manifold $S$ contains two Lagrangian tori $T_1, T_2$ such that $\pi_1(S-(T_1\cup T_2))$ is generated by  the two commuting elements $s_1, t_1$. The Mayer-Vietoris sequence computes $H_1(S-(T_1\cup T_2))=\ZZ^2$ so that  $\pi_1(S-(T_1\cup T_2))=\ZZ s_1 \oplus \ZZ t_1.$

  The  meridians and Lagrangian push offs of $T_1$ and $T_2$ are given by
$\mu_{T_1}=1, m_{T_3}=s_1, \ell_{T_1}=s_1^{-1}$ and $\mu_{T_2}=1, m_{T_2}=t_1, \ell_{T_2}=s_1$.
Thus $(S,T_1, T_2)$ is a telescoping triple.    It is shown to be minimal in the proof of \cite[Theorem 8]{BK4}.  The existence of a square $-1$ torus follows exactly as in the proof of Theorem \ref{manC}, since the manifold $S$ is obtained by Luttinger surgeries on the symplectic sum of $(T^2\times T^2)\#2\bCP^2$ and $(T^2\times S^2)\#4\bCP^2$ along a genus 2 surface.

Relabel $S$ as $D$.
 \end{proof}

 The symplectic 4-manifold $X_{1,7}$ obtained from $D$ by $+1$ Luttinger surgeries on $T_1$ and $T_2$ as in Proposition \ref{surger} is   minimal and homeomorphic to $\CP^2\#7\bCP^2$ (\cite{BK4}).
More generally, the following proposition is true.

\begin{prop}\label{minlut} Let $X$ be one of the manifolds $B,B_g, C,D$ and $T_1, T_2$ the corresponding Lagrangian tori  as described in Theorems \ref{manB}, \ref{manC},
\ref{BTypeII}, with Lagrangian push offs $m_{T_i}$ and $\ell_{T_i}$ (and trivial meridians).

Then the symplectic 4-manifolds obtained from $\pm 1$ Luttinger surgery on one or both of $T_1,T_2$ along $m_{T_i}$ or $\ell_{T_i}$  are all minimal.
\end{prop}

We omit the proof, which is based on Usher's theorem and a repeated use of \cite[Lemma 2]{BK4}. The reader may look at the proofs of Theorems 8, 10, and 13 of \cite{BK4}.

 \bigskip

 Since our emphasis in this article is on 4-manifolds with odd intersection forms, we recall the following theorem due to Gompf \cite{Gompf2}.

 \begin{thm}\label{E(k)} The symplectic manifold  $E'(k)= E(k)_{2,3}$ obtained from the elliptic surface $E(k)$ by performing two log transforms of order $2$ and $3$ is    simply connected, and minimal. It
 has Euler characteristic $e(E'(k))=12k$, signature $\sigma(E'(k))=-8k$, and an odd intersection form.  \qed \end{thm}

 \section{Minimal symplectic 4-manifolds with signature $-3$ and Euler characteristic greater than 14}\label{sig3}

 The most complicated  examples we construct are simply connected minimal symplectic 4-manifolds with signature $-3$. Putting these in the context of telescoping triples is more trouble than constructing them directly. Moreover,  with the exception of the $\sigma=-3$ manifolds,  our inductive scheme for filling out the entire geography for $\sigma<-1$ only requires at most one copy of the manifold $A$ of Lemma \ref{manA}. Hence in this section we  prove the following theorem.

\begin{thm} \label{sig-3} For  each integer $k\ge 2$, there exists a simply connected minimal symplectic 4-manifold $X_{1+2k, 4+2k}$ with $e(X_{1+2k, 4+2k})= 7+4k$ and $\sigma(X_{1+2k, 4+2k})=-3$.
\end{thm}
 The construction of signature $-3$ 4-manifolds for $e=7+8g$ is easier than  for $e=11+8g$.  Roughly speaking, to produce a 4-manifold  with $e=7+8g$, we take the symplectic sum along a genus 3 surface of the 4-manifold $A$ of  Lemma \ref{manA} with $F\times G$, where $F$ is a genus 3 surface and $G$ is a genus $g$ surface, and perform Luttinger surgery on the Lagrangian tori in $G$. To produce a 4-manifold with $e=11+8g$ requires producing a substitute $A'$ for $A$ which has signature $-3$ and $e=11$, and which satisfies the conclusions of Lemma \ref{manA}.  To do this, we take the symplectic sum of $A$ with the product $F\times G$ of two genus 2 surfaces along a symplectic torus.

 \begin{lem}\label{gen2xgen2} There exists a minimal symplectic 4-manifold $Z$ with $e(Z)= 4$ and $\sigma(Z)=0$ which contains  eight homologically essential Lagrangian tori $S_1,\cdots, S_8$ (in fact each $S_i$ has a geometrically dual torus $S_i^d$ so that all other intersections are zero) so that $\pi_1(Z-\cup_i S_i)$ is generated by   $x_1, y_1,x_2,y_2$ and $ a_1,b_1,a_2,b_2$, and so that the meridians and Lagrangian push offs are given by
 \begin{itemize}
\item $S_1: \mu_1= [b_1^{-1}, y_1^{-1}], m_1=x_1, \ell_1=a_1$,
\item $S_2: \mu_2=[x_1^{-1}, b_1], m_2=y_1, \ell_2=b_1  a_1  b_1^{-1},$
\item $S_3: \mu_3= [b_2^{-1},y_1^{-1}], m_3=x_1, \ell_3=a_2$,
\item $S_4: \mu_4= [x_1^{-1}, b_2], m_4=y_1, \ell_4=b_2a_2b_2^{-1},$
\item $S_5: \mu_5= [b_1a_1^{-1}b_1^{-1}, y_2^{-1}], m_5=x_2, \ell_5=b_1^{-1},$
\item $S_6: \mu_6=[x_2^{-1}, b_1a_1b_1^{-1}], m_6=y_2, \ell_6=b_1 a_1 b_1^{-1} a_1^{-1}b_1^{-1},$
\item $S_7: \mu_7= [b_2a_2^{-1}b_2^{-1},y_2^{-1}], m_7=x_2, \ell_7=b_2^{-1}$,
\item $S_8: \mu_8= [x_2^{-1}, b_2a_2b_2^{-1}], m_8=y_2, \ell_8=b_2 a_2 b_2^{-1} a_2^{-1}b_2^{-1}$.
\end{itemize}
 \end{lem}

\begin{proof}  Proposition 7 of  \cite{BK4}  (see also the construction of the manifold $P$ in \cite{BK2}) computes the
fundamental group of the complement of four Lagrangian tori $S_1, S_2, S_3, S_4$   in the product
$F_1\times G$ of a punctured torus $F_1$ with a genus 2 surface $G$.   This group is generated by
loops $x_1,y_1,a_1,b_1,a_2,b_2$  (called $\tx,\ty, \ta_1,\tb_1,\ta_2,\tb_2$ there) where $a_1, b_1, a_2, b_2$ are a standard generating set for $\pi_1(G)$, and $x_1,y_1$ are a standard generating set for $\pi_1(F_1)$ based at a   point $h$ on the boundary.  In particular, the boundary of $F_1\times G$ is
$\partial F_1\times \{k\}\cup \{h\}\times G$, and the copy $\{h\}\times G$ carries the loops $a_1,b_1,a_2,b_2$.

We take two copies of this manifold, calling the second $F_2\times G$, its tori $S_5,S_6,S_7,S_8$, and its generators   $x_2,y_2,a'_1,b'_1,a'_2,b'_2$.

Glue the two copies together using a diffeomorphism of their boundary of the form
$Id\times \phi:\partial H_2\times G\to \partial H_1\times G$, where
$\phi:G\to G$ is the base point preserving diffeomorphism inducing the map
$$(a_1',b_1',a_2',b_2')\mapsto (b_1^{-1}, b_1a_1b_1^{-1},  b_2^{-1}, b_2a_2b_2^{-1})$$
(a composite of six Dehn twists: see \cite[Lemma 9]{BK4}).

\medskip

The resulting manifold $Z$ can also be described as the symplectic sum of two copies of a product of a genus 1 and  genus 2 surface. Thus the result is symplectic and the 8 tori are Lagrangian. The tori $S_1,S_2,S_3,S_4$ in $H_1\times G$ have  geometrically dual tori $S_1^d,S_2^d, S_3^d, S_4^d$  which form a direct sum (geometrically) of four hyperbolic pairs, and similarly for $S_5,S_6,S_7,S_8$.    Clearly $e(Z)=4$ and $\sigma(Z)=0$.

 Applying the Seifert-Van Kampen theorem to the formulae of Proposition 7 of  \cite{BK4}     finishes the fundamental group assertions.

Since the diffeomorphism $ Id\times \phi:\partial H_2\times G \to  \partial H_1\times G$ extends to
$H_2\times G\to H_1\times G$, the manifold $Z$ is nothing but the product of two genus 2 surfaces.  In particular $Z$ is minimal.
 \end{proof}

 \medskip
 Let $Y$ be the symplectic 4-manifold obtained from the manifold $Z$ of Lemma \ref{gen2xgen2}  by performing the following seven Luttinger surgeries on $S_1,\cdots, S_7.$

\begin{enumerate}
\item $S_1: +1$ surgery along $m_1$.
\item $S_2: +1$ surgery along $\ell_2$.
\item $S_3: +1$ surgery along $\ell_3$.
\item $S_4: +1$ surgery along $m_4$.
\item  $S_5: +1$ surgery along $\ell_5$.
\item $S_6: +1$ surgery along $m_6$.
\item $S_7: +1$ surgery along $m_7$.
\end{enumerate}
Since the torus $S_8$ has not been surgered,   it remains as a Lagrangian  torus in $Y$. Since $S_8$ is homologically essential, the symplectic form can be perturbed so that $S_8$ becomes symplectic. The  symplectic 4-manifold $Y$ is minimal, since it is a symplectic sum of manifolds with contractible universal cover (see \cite[Lemma 2]{BK4}).

\bigskip

Let $(B, T_1, T_2)$ be the telescoping triple   of Theorem \ref{manB}, with $B$ containing the genus 2 symplectic surface $F$ and geometrically dual $-1$ torus $H_1$.  Perform $+1 $ Luttinger surgery on $T_2$ along $m_{T_2}$ to kill $t_1$, yielding a minimal  (Proposition \ref{minlut})  symplectic 4-manifold $\hat{B}$. Note that $\hat{B}$ still contains the three  surfaces $T_1, F, H_1$ and $\pi_1(\hat{B}-(T_1\cup F\cup H_1))=\ZZ t_2$.   The torus $T_1$ is disjoint from the geometrically dual symplectic surfaces $F$ and $H_1$, and its Lagrangian push offs are $m_{T_1}=1$ and $\ell_{T_1}=t_2$, by Theorem \ref{manB}.

\medskip

\begin{lem} \label{manB'} The symplectic sum
$X_{3,5}= \hat{B}\#_{T_1, S_8} Y$  is simply connected, minimal,  contains a symplectic genus 2 surface of square 0 and a geometrically dual symplectic torus of square $ -1$.  Moreover, $e(X_{3,5})=10$ and $\sigma(X_{3,5})=-2$, so that $X_{3,5}$ is homeomorphic but not diffeomorphic to $3\CP^2\# 5\bCP^2$.

\end{lem}

 \begin{proof}  We refer to the notation in the statement of Theorem \ref{manB}. The fundamental group of $X_{3,5}-F$ is generated by $  t_2, x_1,y_1,x_2,y_2,a_1,b_1,a_2,b_2$ by the Seifert-Van Kampen theorem (recall that $t_1$ is killed by Luttinger surgery on $T_2$).

 Since the meridian of $T_1$ in $\pi_1(B-(F\cup T_1\cup T_2))$ is trivial, $\mu_8$ is trivial in $\pi_1(X_{3,5}-F)$.   Choose the gluing map $S_8\to T_1$ so that $\ell_8$ is killed and $m_8$ is sent to $t_2$ (i.e. $m_8\mapsto \ell_{T_1}, \ell_8\mapsto m_{T_1}^{-1}$).

Since $\ell_8$ is a conjugate of $b_2^{-1}$ and $m_8=y_2$,  it follows that $b_2=1$ and $y_2=t_2$.  This implies that $\mu_3$ and $\mu_4$ are trivial, and hence the third and fourth Luttinger surgeries listed above show that $a_2=1$ and $y_1=1$. Thus $\mu_1$ and $ \mu_7$ are killed. The first and seventh Luttinger surgeries now show that $x_1=1$ and $x_2=1$. Continuing, we see that $\mu_2$ and $\mu_6$ are killed so that the corresponding surgeries give $a_1=1$ and $y_2=1$. This implies $\mu_5=1$ and so $b_1=1$.  Hence $\pi_1(X_{3,5})=1$.

That $X_{3,5}$ is minimal follows from Usher's theorem.  The genus 2 surface $F$ and torus $H_1$ in $B$ survive to give the required surfaces in $X_{3,5}$.
 \end{proof}

 \medskip

 \medskip

 \noindent{\em Proof of Theorem \ref{sig-3}.} We define two minimal simply connected symplectic 4-manifolds: let
  $X_{-}= X_{3,5}$  and let $X_{+}=X_{1,3}$  (thus $X_{+}$ is obtained from the manifold $\hat{B}$ defined above by performing $+1$ Luttinger surgery on $T_1$ along $\ell_{T_1}$; see Corollary \ref{1-3}).
Then $X_{-}$ and $X_{+}$ each contain  a symplectic genus 3 surface $F_3$ of square 1 obtained by resolving the union $H_1\cup F$.  Moreover, $e(X_{-})=10$, $\sigma(X_{-})=-2$, $e(X_{+})=6$,  and $\sigma(X_{+})=-2$.

  Blow up $X_\pm$ once at a point on $F_3$ and take the proper transform. Call the result $\tilde{X}_\pm$ and denote by  $\tilde{F}_3 $ the proper transform of $F_3$.
  Thus $\tF_3$ is a genus 3, square zero symplectic surface with simply connected complement, which meets every $-1$ sphere in $\tilde X_\pm$ since $X_\pm$ is minimal.

 We now mimic  the proof of Corollary \ref{cor6}. Take the product $F_3\times G$ of a genus 3 surface with a genus $g$ surface. Perform Luttinger surgeries on the $2g$ disjoint Lagrangian tori $Y_1\times A_j$ and $Y_2\times B_j$  along the curves $\ell_{Y_1\times A_j}=a_j$  and $\ell_{Y_2\times B_j}=b_j$ to obtain a manifold $Z_g$.

 Then by Theorem  \ref{thm:technical_lem}  the fundamental group of $Z_g$ is generated by  the $6+2g$ loops $\tx_1,\ty_1,$ $\tx_2,$ $\ty_2,$ $\tx_3, $ $ \ty_3$,  $\ta_1,\tb_1,\cdots, \ta_g,\tb_g$, and the relations
 $$[\tx_1,\tb_j]=\ta_j, [\tx_2,\ta_j]=\tb_j$$ hold in $\pi_1(Z_g)$. Moreover, the standard symplectic generators for $\pi_1(F)$ are sent to $\tx_1,$ $\ty_1,$ $\tx_2,$ $\ty_2,$ $ \tx_3,\ty_3$ in $\pi_1(Z_g)$.

 Since $\pi_1(\tilde X_\pm-\tF_3)=1$, the fundamental group of the symplectic sum
 $$ Q_{\pm, g}=\tilde X_\pm\#_{\tF_3,F}Z_g$$
 is trivial. Indeed, the $\tx_i,\ty_i$ are killed by taking the symplectic sum, and the relations coming from the Luttinger surgeries show the $a_i$ and $b_i$ are killed also.

 Now $Q_{\pm,g}$ is minimal provided $g\ge 1$ by Usher's theorem since $\tilde X_\pm$ is relatively minimal by Li's theorem, Theorem \ref{Li}.

 One computes:
 $$e(Q_{-,g})=e(\tilde X_{-})+e(Z_g)+8= 11+ 8g-8+8=11+ 8g, $$
 $$\sigma(Q_{-,g})=
 \sigma(\tilde X_{-})+\sigma(Z_g)=  -3$$
 and
 $$e(Q_{+,g})=e(\tilde X_{+})+e(Z_g)+8= 7+ 8g-8+8=7+ 8g, $$
 $$\sigma(Q_{+,g})=
 \sigma(\tilde X_{+})+\sigma(Z_g)= -3 $$

 Thus we set $X_{1+2k, 4+2k}= Q_{+, k/2}$ if $k$ is even and $X_{1+2k, 4+2k}= Q_{-, (k-1)/2}$ if $k$ is odd.

 This completes the proof of Theorem \ref{sig-3}.
 \qed

 \medskip

 \noindent{\bf Remark.}    In the construction of the manifold   $X_{5,8}$,   the first step (see the paragraph preceding Lemma \ref{manB'}) involves Luttinger surgery on the torus $T_2$ to kill $t_1$.  If one constructs the  manifold  $P_{5,8}$   by by the same construction as for $X_{5,8}$ except by {\em not} performing this surgery, then $P_{5,8}$ is  a minimal symplectic 4-manifold with
 $\pi_1(P_{5,8})=\ZZ t_1$ containing an essential Lagrangian (or if desired symplectic) torus  $T=T_2$ such that the inclusion map $\pi_1(T)\to \pi_1(P_{5,8})$ is a surjection and the inclusion map $\pi_1(P_{5,8}-T)\to \pi_1(P_{5,8})$ an isomorphism. Moreover $e(P_{5,8})=15$ and $\sigma(P_{5,8})=-3$.

More generally, for any $k\ge 2$ the same construction yields a minimal symplectic 4-manfold $P_{1+2k, 4+2k}$ containing a Lagrangian or symplectic torus $T$ with these properties and  such that
$e(P_{1+2k, 4+2k})=7+4k$, $\sigma(P_{1+2k, 4+2k})=-3$.

\section{Small examples with odd signature}\label{small}

In this section, we remind the reader of some known examples of small manifolds with odd signature, and construct a few  new ones.

 Kotschick showed   in \cite{Kot}  that the Barlow surface  is smoothly irreducible  and hence it is a minimal  symplectic 4-manifold homeomorphic to $\CP^2\#8\bCP^2$. This manifold realizes the pair
 $e=11, \sigma=-7$.

In \cite{Gompf}, Gompf constructs small minimal symplectic 4-manifolds   which contain appropriate tori. For example, the manifold Gompf calls $S_{1,1}$ is minimal, has $e = 23$ and $\sigma= -15$,  and contains a symplectic torus of square zero with simply connected complement (\cite[Lemma 5.5]{Gompf}). The minimality  $S_{1,1}$ was proved by Stipsicz \cite{S2}.

Gompf also constructs other minimal   symplectic 4-manifolds: the manifold $R_{2,1}$ has    $e= 21$ and $\sigma= -13$ and $R_{2,2}$ has $e=19$ and $\sigma=-11$.   The minimality of $ R_{2,1}$ was proved by J. Park \cite{JP2}, and $R_{2,2}$ was proven to be minimal by  Szab\'o \cite{Sz2}.

In   \cite{SS}, Stipsicz and Szabo construct a minimal   symplectic 4-manifold homeomorphic to  $ \CP^2\#6\bCP^2$, realizing $ e=9,\sigma=-5$.

In \cite{DP}, the fifth author constructs a minimal simply connected symplectic 4-manifold  homeomorphic to $3\CP^2\# 12\bCP^2$, hence with $e=17$ and $\sigma=-9$,  containing a symplectic torus $T_{2,4}$ with simply connected complement. This manifold is called $X_{12}$ in that article, we will use the notation $X_{3,12}$ here to avoid confusion.

\bigskip

We produce   a few more small examples.

\begin{prop}\label{5-10} There exists a minimal simply connected symplectic 4-manifold $X_{5,10}$ homeomorphic to $5\CP^2\# 10\bCP^2$, hence with $e=17$ and $\sigma=-5$.
\end{prop}
\begin{proof} The manifold $X_{1,3}$ of Corollary \ref{1-3} contains a symplectic genus 2 surface $F$ of square zero, and a geometrically dual symplectic torus $H_1$ with square $-1$. Symplectically resolve $F\cup H_1$ to produce a square $1$ symplectic genus 3 surface $F_3\subset X_{1,3}$.

Blow up $X_{1,3}$ at a point on $F_3$ to obtain $\tilde{X}_{1,3}$ and take $\tilde{F}_3$ to be the proper transform of $F_3$. Then  $\tilde{F}_3$ is a square zero symplectic surface that meets every $-1$ sphere in $\tilde{X}_{1,3}$ by Li's theorem (Theorem \ref{Li}).   Moreover, since $X_{1,3}$ is simply connected and $\tilde{F}_3$ meets the exceptional sphere, $\tilde{X}_{1,3}-\tilde{F}_3$ is simply connected.

Take $Y=T\times F_2$, the product of a torus with a genus 2 surface. Then $Y$ contains the geometrically dual symplectic surfaces $T\times \{p\}$ and $\{q\}\times F_2$. Symplectically resolve their union to obtain a genus 3, square 2 symplectic  surface $F_3'\subset Y$.  Note that the homomorphism induced by  inclusion $\pi_1(F_3')\to \pi_1(Y)$ is surjective. Blow up $Y$ twice at points on $F_3'$ to obtain $\tilde{Y}$ and the proper transform $\tilde{F}_3'$, a square zero genus 3 symplectic surface.

Then the symplectic sum
$$X_{5,10}= \tilde{X}_{1,3}\#_{\tilde{F}_3,\tilde{F}'_3}\tilde{Y}$$
is simply connected. It is minimal by Usher's theorem.

Its characteristic classes are computed
$$e(X_{5,10})=e( \tilde{X}_{1,3})+e(\tilde{Y})+8=7+2+8=17$$
and
$$\sigma(X_{5,10})=\sigma( \tilde{X}_{1,3})+\sigma(\tilde{Y})=-3-2=-5.$$
The proposition follows.
\end{proof}

\begin{prop}\label{5-12} There exists a minimal simply connected symplectic 4-manifold $X_{5,12}$
homeomorphic to $5\CP^2\# 12\bCP^2$, hence with $e=19$ and $\sigma=-7$.

\end{prop}
\begin{proof} The proof is very similar to the proof of Proposition \ref{5-10}. Construct $\tilde{X}_{1,3}$ and $\tilde{F}_3$ as in that proof.

Take  $Z=T\times T$, the product of two tori.  Pick three distinct points $p_1,p_2, q$ in $T$.  Then $Z$ contains the  three symplectic surfaces $T\times \{p_1\}$,  $T\times \{p_2\}$ and $\{q\}\times T$.  Symplectically resolve their union to obtain a genus 3, square 4 symplectic  surface $F_3'\subset Z$.  Note that the homomorphism induced by  inclusion $\pi_1(F_3')\to \pi_1(Z)$ is surjective. Blow up $Z$ four times at points on $F_3'$ to obtain $\tilde{Z}$ and the proper transform $\tilde{F}_3'$, a square zero genus 3 symplectic surface.

Then the symplectic sum
$$X_{5,12}= \tilde{X}_{1,3}\#_{\tilde{F}_3,\tilde{F}'_3}\tilde{Z}$$
is simply connected. It is minimal by Usher's theorem.

Its characteristic classes are computed
$$e(X_{5,12})=e( \tilde{X}_{1,3})+e(\tilde{Z})+8=7+4+8=19$$
and
$$\sigma(X_{5,12})=\sigma( \tilde{X}_{1,3})+\sigma(\tilde{Z})=-3-4=-7.$$
The proposition follows.

\end{proof}

\begin{prop}\label{5-14} There exists a minimal simply connected symplectic 4-manifold $X_{5,14}$
homeomorphic to $5\CP^2\# 14\bCP^2$, hence with $e=21$ and $\sigma=-9$.

\end{prop}
\begin{proof} The proof is very similar to the proof of Proposition \ref{5-12}. Construct $\tilde{X}_{1,3}$ and $\tilde{F}_3$ as in that proof.

Take  $Z=T\times S^2$, the product of a torus and a sphere.  Pick three distinct points $p_1,p_2, p_3 $ in $S^2$ and $q\in T$.  Then $Z$ contains the  four symplectic surfaces $T\times \{p_1\}$,  $T\times \{p_2\},
T\times \{p_2\}$ and $\{q\}\times S^2$.  Symplectically resolve their union to obtain a genus 3, square 6 symplectic  surface $F_3'\subset Z$.  Note that the homomorphism induced by  inclusion $\pi_1(F_3')\to \pi_1(Z)$ is surjective. Blow up $Z$ six times at points on $F_3'$ to obtain $\tilde{Z}$ and the proper transform $\tilde{F}_3'$, a square zero genus 3 symplectic surface.

Then the symplectic sum
$$X_{5,14}= \tilde{X}_{1,3}\#_{\tilde{F}_3,\tilde{F}'_3}\tilde{Z}$$
is simply connected. It is minimal by Usher's theorem.

Its characteristic classes are computed
$$e(X_{5,14})=e( \tilde{X}_{1,3})+e(\tilde{Z})+8=7+6+8=21$$
and
$$\sigma(X_{5,14})=\sigma( \tilde{X}_{1,3})+\sigma(\tilde{Z})=-3-6=-9.$$
The proposition follows.
\end{proof}

\section{The main theorem}

In this section we prove the first  theorem stated in the introduction.  We begin   with an arithmetic lemma.  The purpose of this lemma is to produce the number of     each of the model manifolds $B, B_g, C, D, E(k)$ needed to construct 4-manifolds with specified  signature and Euler characterstic. The proof includes an algorithm for finding these numbers.

\medskip

\begin{lem} \label{arith} Given any pair of non-negative  integers $(m,n)$ such that
$$0\leq m\leq 4n-1$$
there exist non-negative integers $b,c,  d,g,$ and $k$ so that
$$m= d+2c+3b+4g \text{ and } n=b+c+d+k+g$$
and so that $b\ge 1$ if $g>0$.
 \end{lem}

\begin{proof} If $m=0$, set $k=n$ and $b=c= d=g=0$.

Assume then that $m>0$. Choose a non-negative integer $\ell$ so that $\frac{m+1}{4}\leq n-\ell\leq m$.  Let
$$s=\max\{z\in\ZZ\ |  \ 3z\leq 4n-4\ell-m-1\} \ \text{ and } \ \Delta= 4n-4\ell-m-1-3s.$$
Then $\Delta=0,1$ or $2$, and $s\ge 0$. Moreover
$$n-\ell-1-s=\tfrac{1}{3}(m-(n-\ell))-\tfrac{2}{3}+\tfrac{\Delta}{3}\ge \tfrac{\Delta}{3}-\tfrac{2}{3}.$$

If $\Delta=0$, then set $b=1,c=0,  d=s, g=n-\ell-s-1$, and $k=\ell$.  Since $\Delta=0$,
and since $g$ is an integer, $g\ge 0$.

If $\Delta=1$ and $s\ge 1$, then
 set $b=1, c=0, d=s-1, g=n-\ell-s-1$, and $k=\ell+1.$ Note that $g\ge -\tfrac{1}{3}$ so that $g\ge 0$.

 If  $\Delta=1$ and $s=0$, then   either $n-\ell-2\ge0$   in which case we set $b=2, c=0, d=0, k=\ell, $ and $g=n-\ell-2$, or else $n-\ell-2=-1$ in which case
we take $b=0, c=1, d=0, k=\ell$, and $g=0$.

If $\Delta=2$ and $s\ge 2$, set $b=1,c=0,d=s-2,g=n-\ell-s-1$ and $k=\ell+2$.

If $\Delta=2$, $ s=1$, and $n-\ell\ge 3$ then set $b=2, c=0, d=0, g=n-\ell-3$ and $k=\ell+1$. If $\Delta=2$, $ s=1$, and $n-\ell< 3$, then  necessarily $n-\ell=2$, and so $(m,n)=(2,\ell+2)$ and we set $b=0,c=1, d=0,k=\ell+1,g=0$.

This leaves the cases when $\Delta=2$ and $s=0$. If $n-\ell \ge 2$, set $b=1,c=1, d=0, g=n-\ell-2, $ and $k=\ell$. Finally, if $n-\ell=1$, then $(m,n)=(1, 1+\ell)$, so we take $b=0,c=0,d=1, g=0$ and $k=\ell$.
\end{proof}

We can now prove our main result.  We state it in terms of $c_1^2 = 2e +3\sigma $ and $\chi_h=\tfrac{1}{4}(e+\sigma)$ because it is simpler to work with these numbers than pairs $(e,\sigma)$ where $e+\sigma =0 \mod 4$.  Note that in this notation, a 4-manifold with $c_1^2=8\chi_h+k $ has signature $k$, so the line $c_1^2=8\chi_h-2$ corresponds to manifolds with signature $-2$.

\begin{thm}\label{geosimp} For any pair $(c, \chi)$ of non-negative integers satisfying
$$0\leq c \leq 8\chi -2 $$  with the   possible exceptions of
$
(c , \chi )=  (5,1)$, $(9,2)$, $(13,2)$,  or $ (13,2)$,
 there exists a  minimal simply connected symplectic 4-manifold $Y=X_{2\chi-1, 10\chi-c-1}$ with odd intersection form and
$$c = c_1^2(Y)  \text{ and } \chi =\chi_h(Y).$$
Hence $Y$   is   homeomorphic but not diffeomorphic  to
$(2\chi-1)\CP^2\#(10\chi-c-1)\bCP^2.$\end{thm}
\begin{proof}

We make extensive use of the manifolds $A,B, B_g, C, D, E'(k)$ of (respectively) Lemma \ref{manA}, Theorem \ref{manB}, Corollary \ref{cor6}, Theorem \ref{manC}, Theorem \ref{BTypeII}, and Theorem \ref{E(k)}. We will also use the sporadic examples of Section \ref{small}.

\medskip

We first realize all  pairs with $c$ even.
Let $(m,n)=(\tfrac{1}{2}c, \chi)$.  Lemma \ref{arith} produces integers $b,c,d,g, $ and $k$ so that $m= d+2c+3b+4g$  and $n=b+c+d+k+g$
and so that $b\ge 1$ if $g>0$.

 Construct the symplectic sum along tori of
\begin{enumerate}
\item $b$ copies of $B$  if $g=0$, or one copy of $B_g$ and $b-1$ copies of $B$ if $g\ge 1$,
\item $c$ copies of $C$, and
\item $d$ copies of $D$.
 \end{enumerate}
More  precisely, each of the manifolds $B,C,D$ contain two essential Lagrangian tori. Construct the
symplectic sum $Z$ of these manifolds  by chaining them together, using Proposition \ref{symp_sum_type_I_or_II} to ensure that at each stage one has a telescoping triple.

Specifically, if $g=0$ take
$$Z= B\#_s  \cdots \#_s B\#_sC\#_s\cdots \#_s C\#_s D\#_s\cdots \#_s D$$ and if $g\ge 1$ take
$$Z=B_g\#_sB\#_s  \cdots \#_s B\#_sC\#_s\cdots \#_s C\#_s D\#_s\cdots \#_s D$$
where $\#_s$ denotes the symplectic sum along the appropriate tori (perturbing the symplectic forms so that they become symplectic)  according to the recipe of Theorem \ref{symp_sum_type_I_or_II},   so that the two unused Lagrangian tori (which we relabel $T_1$ and $T_2$) make $(Z,T_1,T_2)$ a telescoping triple.

If $k=0$, then perform $+1$ Luttinger surgery on $T_1$ and $T_2$ to obtain a simply connected (according to Proposition \ref{surger})   symplectic 4-manifold $Y$.

If $k\ge1$ and one of $b,c,d$ is positive,   perform $+1$ Luttinger surgery on $T_2$ in $Z$ and take the symplectic sum of the result with the elliptic surface $E(k)$  along $T_1$ to obtain the manifold $Y$. Since $E(k)-T$ is simply connected, so is $Y$, by the same reasoning as in the proof of Proposition \ref{surger}. Since $B$, $C$, and $D$ contain $+1$ tori disjoint from the Lagrangian tori $T_1, T_2$, the manifold $Y$ has odd intersection form.

If $k\ge 1$ and $b,c,d$ are zero, take $Y=E'(k)$ (see Theorem \ref{E(k)}) which has an odd intersection form.

Thus $Y$ is a simply connected symplectic manifold realizing the pair $(c, \chi)$. Since each of the manifolds $B,B_g, C,$ and $D  $   contain a surface of odd square which misses the tori used in forming the symplectic sums, and since $E'(k)$ has an odd intersection form, it follows that $Y$ has an odd intersection form.

 Since the 4-manifold $Y$ has indefinite, odd intersection form, Freedman's theorem \cite{Freedman} implies that  $Y$ is homeomorphic to an appropriate connected sum of $\CP^2$s and $\bCP^2$s.

\bigskip

Now we turn to the case when $c$ is odd. Suppose first that $1\leq c\leq 8\chi -17$.
Let $(c',\chi')= (c-1, \chi-2)$. Thus $0\leq c'\leq 8\chi'-2$, and $c'$ is even. Construct the manifold   $Z$ corresponding to  the pair $(c',\chi')$ and either perform $+1$ Luttinger surgery on $T_1$ or take the symplectic sum with $E(k)$ if $k\ge 1$.  But rather than performing $+1$ Luttinger surgery on $T_2$ as we did above, perturb the symplectic form to make $T_2$ symplectic, and then take the symplectic sum with  Gompf's  manifold $S_{1,1}$ (see Section \ref{small}) along the symplectic torus in $S_{1,1}$ with simply connected complement. Since $S_{1,1}$ has $c_1^2=1$ and $\chi_h=2$ the resulting symplectic manifold $Y$ has
$(c_1^2,\chi_h)=(c,\chi)$.

Next suppose that $c$ is odd and $7\leq c\leq 8\chi-11$. Set $(c',\chi')=(c-7,\chi-2)$. Thus $0\leq c'\leq 8\chi'-2$ and $c'$ is even. Construct the manifold   $Z$ corresponding to  the pair $(c',\chi')$.
We repeat the argument of the previous paragraph, replacing Gompf's manifold $S_{1,1}$ with the manifold  $X_{3,12}$   of Section \ref{small}. Take the symplectic sum of $Z$ with $X_{3,12}$ along $T_2$ and $T$.  Since $ c_1^2(X_{3,12})=7$ and $\chi_h(X_{3,12})=2$,     the resulting manifold $Y$ realizes  the pair $(c,\chi)$.

To  realize  all pairs $(c,\chi)$ with $c$ odd and $21\leq c\leq
8\chi-5$, repeat the argument once more, this time using the
manifold $P_{5,8}$ described in the remark at the end of the proof
of Theorem \ref{sig-3}, which has $c_1^2=21$ and $\chi_h=3$. A bit
of care must be taken to ensure that the  result is simply connected
since $\pi_1(P_{5,8})=\ZZ$. This is accomplished by making sure that
the generator of $\pi_1(T)$ sent to the generator of
$\pi_1(P_{5,8}-T)$ is identified with an element in the kernel of
$\pi_1(T_2)\to \pi_1(Z-T_2)$ when forming the sympletic sum $Y=Z\#_s P_{5,8}$.

The manifold  $Y=X_{1+2k, 4+2k}$ of Theorem \ref{sig-3} provides an example  realizing $(c,\chi)=(5+8k, 1+k)$ for any $k\ge 2$, i.e. $21\leq c=8\chi-3$.

Since $c_1^2\equiv \sigma\pmod{2}$, and simply connected 4-manifolds with odd signature have an odd intersection form, it follows that the manifolds constructed for $c$ odd also have an odd intersection form.

 \medskip

It remains to show that $Y$ is minimal.   Since $E'(k)$  is minimal, we assume that $c_1^2>0$.    By Proposition \ref{minlut}, the 4-manifold obtained by performing one or two  $\pm 1$ Luttinger surgeries on  $T_1$ or $T_2$ along $\ell_{T_i}$ or $m_{T_1}$ in $B, C,$ or $D$ is minimal.  The $E(k)$ are minimal for $k\ge 2$. Although $E(1)$ is not minimal, every $-1$ sphere intersects the elliptic torus. Thus   $Y$ is the symplectic sum of minimal (or, if $k=1$, relatively minimal) symplectic 4-manifolds and therefore is minimal by Usher's theorem.

It is easy to check that the only pairs $(c, \chi)$ with $0\leq c\leq 8\chi-2$ which are omitted by these cases are

$$\begin{array}{l}
(1,1), (3,1), (5,1),\\[.2cm]
(1,2), (3,2), (5,2), (7,2), (9,2), (11,2), (13,2),\\[.2cm]
(15,3), (17,3), (19,3).\\[.2cm]
\end{array}$$

The examples listed in Section \ref{small} realize most of these pairs. The only ones left unrealized are $(5,1)$, $(9,2)$, $(13,2)$, and $(13,2).$

 \end{proof}

The four unrealized pairs do correspond to (non-minimal) symplectic 4-manifolds; e.g. blow ups of  $X_{1,3}$ or $X_{3,5}$.

 It is conjectured that the irreducible smooth 4-manifold homeomorphic to $3\CP^2\#10\bCP^2$
constructed in \cite{DP} and the irreducible smooth 4-manifold homeomorphic to $3\CP^2\#8\bCP^2$
constructed in \cite{JP3} are symplectic (and hence minimal): their Seiberg-Witten invariants have the right form to be the invvariants of a symplectic manifold.

 There exist small simply connected minimal symplectic 4-manifolds with non-negative signature (e.g. $\CP^2, S^2\times S^2$).  To date, no small examples are known that contain a suitable Lagrangian torus  for which we can extend the construction of Theorem \ref{geosimp}.  Some moderately large examples are known and we will briefly explore the consequences for the geography problem below.

 \bigskip

\noindent{\bf Remark.} Each of the manifolds   constructed in Theorem
\ref{geosimp}, with the possible exception of those corresponding
$c_1^2=0$ and some of the small manifolds with $c$ odd, contain
nullhomologous tori suitable for altering the differentiable
structure  as explained in \cite{FSP}, using \cite{MMS} to compute
the change in Seiberg-Witten invariants. Those with $c_1^2=0$  are $E'(k)$, for which the methods of \cite{GM, FS6, FS8} show how to alter the differentiable structure.  Hence  the manifolds of Theorem
\ref{geosimp} admit infinitely many smooth structures.

\bigskip

The proof of Lemma \ref{arith}  provides a specific algorithm for constructing  simply connected minimal 4-manifolds with desired characteristic numbers, using the model manifolds $A$,$B$, $B_g$, $C$,$D,$ and $E(k)$.

For example, to construct a minimal symplectic manifold homeomorphic but not diffeomorphic to $3\CP^2\#17 \bCP^2$, one sees that such a manifold would have $(c_1^2,\chi_h)=(2,2)$. This corresponds to $(m,n)=(1,2)$ in Lemma \ref{arith}. In the notation of Lemma \ref{arith}, we see that in this case $\ell=1, s=0$ and $\Delta=2$, so that $b=0,c=0,d=1,g=0$, and $k=1$. Thus the desired manifold is obtained by taking the symplectic sum
$$D\#_s E'(1)$$
and performing $+1$ Luttinger surgery on the remaining Lagrangian torus in $D$.

As another example, we construct a minimal symplectic manifold homeomorphic but not diffeomorphic to $21\CP^2\#31 \bCP^2$, i.e.  $\chi_h=11$ and $c_1^2= 78$. Thus $(m,n)=(39, 11)$. The proof of Lemma \ref{arith} provides $\ell=0$, $s=1$ and $\Delta=1$, and so $b=1, c=0, s=1, g=9$ and $k=1$.   Thus the desired manifold is obtained by taking the symplectic sum
$$B_9\#_s E'(1)$$ and performing $+1$ Luttinger surgery on the remaining Lagrangian torus.

The integers produced by the algorithm of  Lemma \ref{arith} are not unique, for example, the choice  $b=2,c=2,d=0,g=8, $ and $k=0$ yields a manifold
$$ B\#_sB_8\#_s C\#_sC.$$
Performing two $+1$ Luttinger surgeries to this manifold yields a (possibly different) minimal symplectic manifold homeomorphic to but not diffeomorphic to $21\CP^2\#31 \bCP^2$.

 \section{Signature greater than $-2$}

 Finding small minimal symplectic 4-manifolds with signature greater than $-2$ poses a special challenge.
 Stipsicz \cite{stip} shows how to produce simply connected minimal symplectic 4-manifolds with positive signature.    The following theorem  provides a method for producing many examples, given one.   It is also useful in studying the geography problem for non-simply connected 4-manifolds.

 To avoid an overly technical statement, we separate the cases of $c$ odd and even, but a more complete statement would have the same  hypotheses on $(c,\chi)$ as in Theorem \ref{geosimp}.

\begin{thm}\label{wedge}  Let $X$ be a symplectic 4-manifold and  suppose that $X$ contains a symplectic torus $T$    such that the   homomorphism  $\pi_1(T)\to \pi_1(X)$ induced by inclusion
is trivial. Then
for any pair $(c , \chi )$ of non-negative integers  satisfying
$$
  0\leq c \leq 8\chi -2   \text{ if $c$ is even,} $$
$$ 1\leq c \leq 8\chi -7   \text{ if $c$ is odd} $$
  there exists a  symplectic  4-manifold $Y$ with $\pi_1(Y)=\pi_1(X)$,
$$c_1^2(Y)=c_1^2(X)+c \text{ and } \chi_h(Y)=\chi_h(X)+\chi.$$
Moreover, if $X$ is minimal (or more generally if $(X,T)$ is relatively minimal)   then the manifold $Y$ is minimal and has an odd, indefinite
intersection form.
\end{thm}

\begin{proof} The argument is the same as the proof of Theorem \ref{geosimp} save for the last step.
Let $Z$ be as in the proof of Theorem \ref{geosimp}.  If $k=0$, then do $+1$ Luttinger surgery on $T_1$ to get a minimal (by Proposition \ref{minlut} and Usher's theorem) manifold $Z_1$ with $\pi_1(Z_1)\cong \ZZ$ containing a symplectic torus $T_2$ (after perturbing the symplectic structure) so that the induced map $\pi_1(T_2)\to \pi_1(Z_1)$ is a split surjection. If $k\ge 1$ then take a fiber sum of $Z$ with $E(k)$ to again get a manifold $Z_1$ with $\pi_1(Z_1)\cong \ZZ$ containing a symplectic  torus $T_2$ so that the induced map $\pi_1(T_2)\to \pi_1(Z_1)$ is a   surjection.

Since the meridian of $T_2$ is nullhomotopic in $Z_1$, the symplectic sum of $Z_1$ and $X$ has fundamental group isomorphic to that of $X$, since the homomorphism $\pi_1(T)\to \pi_1(X)$ is trivial.

Minimality follows as in the proof of Theorem \ref{geosimp} using Usher's theorem.
Since $c_1^2$ and $\chi_h$ are both additive with respect to symplectic sums along tori, the result follows.
\end{proof}

 \medskip

  We will require one more useful fact about  $B$ and $X_{1,3}$ not mentioned in Theorem \ref{manB} or Corollary \ref{1-3}, namely, the existence of a genus 2 square zero symplectic surface $G$ geometrically dual to $F$.  We indicate how to find $G$: $X_{1,3}$ is obtained by Luttinger surgery on 8 Lagrangian tori in
  the symplectic sum of the twice blown up 4-torus $(T^2\times T^2)\#2\bCP^2$ and the product $T\times F_2$ of a torus and a genus 2 surface.

  This symplectic sum is taken along the genus 2 surface in $(T^2\times T^2)\#2\bCP^2$ obtained by resolving $T^2\times\{p\}\cup \{q\}\times T^2$ and blowing up twice (for definiteness at points on  $T^2\times\{p\}$).  In  $T\times F_2$ one takes the surface $\{x\}\times F_2$.

  The square $-1$ torus $H_1$ of Theorem \ref{manB} and Corollary \ref{1-3} was obtained by taking the torus of the form $T\times \{z\}$ which matches up with one of the exceptional spheres in the symplectic sum.  To find the surface $G$, take another nearby torus of the form $T\times \{z'\}$ in $T\times F_2$ and match it up with a torus of the form $\{q'\}\times T^2$.   This is the required surface $G$. (The surface $F$ is a parallel copy of $\{x\}\times F_2$).

 \begin{thm} \label{positive}
 For all integers $k\ge 45$, there exists a simply connected minimal symplectic 4-manifold $X_{2k+1, 2k+1}$ with Euler characteristic $e=4k+4$ and signature $\sigma=0$.

 For all integers $k\ge 49$ there exists a simply connected minimal symplectic 4-manifold $X_{2k-1,2k}$
with Euler characteristic $e=4k+1$ and signature $\sigma=-1$.

  \end{thm}
 \begin{proof} Start with  the telescoping triple $(B, T_1, T_2)$ of Theorem \ref{manB}. It contains a genus 2 square zero symplectic surface $F$ and a geometrically dual square zero symplectic genus 2 surface $G$. The union $F\cup G$ is disjoint from $T_1\cup T_2$.

 Perform $+1$ Luttinger surgery on $T_1$ along $\ell_{T_1}$ to kill  $t_2$.  Call the result $R$.
 Perturb the symplectic form on $R$ slightly so that $T_2$ becomes symplectic. Note that $\pi_1(R-T_2)=\pi_1 (R)=\ZZ t_1$,$\pi_1(T_2)\to \pi_1(R)$ is surjective, and $R$ is minimal  (Proposition \ref{minlut}).

 In \cite[Theorem 18]{BK4}, a minimal symplectic 4-manifold $\tilde X_{3,5}$ homeomorphic to $3\CP^2\# 5\bCP^2$ and containing a pair of symplectic tori $T_3, T_4$ with simply connected complement is constructed.  The symplectic sum $Q=R\#_{T_2,T_3} \tilde X_{3,5}$  is minimal by Usher's theorem. Moreover, $Q$ is simply connected, since   $T_2\subset R$ induces a surjection on fundamental  groups. The surfaces $F$ and $G$ persist as square zero, symplectic geometrically dual surfaces. Since $e(Q)=16$ and $\sigma(Q)=-4$, $Q$ is neither rational nor ruled.

 Notice that the symplectic torus $T_4$ in $Q$ has simply connected complement.

 In $Q$, take $8$ parallel copies of the genus 2 surface $F $
  and one copy of $G$ and symplectically resolve to obtain a genus 18 surface $\Sigma \subset Q$ of square $16$.  Blow up $Q$  16 times, yielding a genus 18 square zero surface $\tilde\Sigma \subset \tilde Q =Q\#16\bCP^2$.    By Li's theorem, every $-1$ sphere in $\tilde{Q}$ intersects $\tilde \Sigma $. Moreover
  $\pi_1(\tilde Q-\tilde{\Sigma})=1$.

 In \cite[Lemma 2.1]{stip}, a Lefschetz fibration $H\to K$ over a surface $K$ of genus 2 is constructed which has $e=75$ and $\sigma=25$. This fibration  admits a symplectic section of square $-1$ and has fiber genus $16$. The 4-manifold $H$ is an algebraic surface, and by the BMY inequality \cite{BPV} is holomorphically minimal. By \cite{HK}, it is also symplectically minimal.  Moreover, $H$ is neither rational nor ruled since it lies on the BMY line.

 Let $\Sigma'\subset H$ denote the symplectic surface obtained by symplectically resolving the union of the fiber and section.   Then $\Sigma'$ has square $1$, and the exact sequence of fundamental groups for a Lefschetz fibration shows that $\pi_1(\Sigma')\to \pi_1(H)$ is surjective.  Blow up $H$ once along $\Sigma'$ and take the proper transform to obtain a square zero, genus 18 surface in $\tilde\Sigma'\subset \tilde H=H\#\bCP^2$  so that $\pi_1(\tilde H - \tilde\Sigma')\to \pi_1(\tilde H)$ is an isomorphism and $\pi_1(\tilde \Sigma')\to \pi_1(\tilde H)$ is surjective.  By Li's theorem, Theorem \ref{Li}, every $-1$ sphere in $\tilde H$ intersects $\tilde\Sigma'$, since $H$ is not rational or ruled.

  Hence the symplectic sum $ S=\tilde Q \#_{\tilde\Sigma, \tilde\Sigma'}\tilde H$ is minimal.  It is simply connected since $\pi_1(\tilde Q-\tilde{\Sigma})=1$ and $\pi_1(\tilde \Sigma')\to \pi_1(\tilde H)$ is surjective.  Moreover, the symplectic torus $T_4\subset S$ has simply connected complement.

Since $S$  is the symplectic  sum along a genus 18 surface,
$$e(S) =e(\tilde Q)+ e(\tilde H)  +4(18-1)=176$$
and
$$\sigma(S)=\sigma(\tilde Q) + \sigma(\tilde H)= 24 -20= 4.$$

Thus $c_1^2(S)=364$ and $\chi_h(S)=45$. It contains the symplectic torus $T_4$ with simply connected complement. Hence Theorem \ref{wedge} establishes the existence of minimal, simply connected  symplectic 4-manifolds
$$X_{89+2\chi, 85+10\chi-c}$$
with $c_1^2=364+c$ and $\chi_h=45+\chi$ for any $(c,\chi)$ satisfying $0\leq c\leq 8\chi - 2$ when $c$ is even.

Taking $c=8\chi-4$ for any $\chi\ge 1$ yields $X_{89+2\chi, 89+2\chi}$, a minimal simply connected symplectic 4-manifold with signature zero. The intersection form is odd since, as one can check from Lemma \ref{arith}, either $\chi=1$ in which case  the model manifold $C$ (with its $-1$ torus) is used in the construction of $X_{91,91}$, or else $\chi>1$, in which case the model manifold $B $ (with its $-1$ torus) is used in the construction of $X_{89+2\chi,89+2\chi}.$

  \medskip

  To get minimal symplectic 4-manifolds with signature $-1$ , consider the symplectic sum
  $$Y= B\#_{T_1, T}P_{1+2k, 4+2k}$$ of the manifold $B$ of Theorem \ref{manB} with the manifold  $P_{1+2k, 4+2k}$ of the remark at the end of Section \ref{sig3} along $T_1$ in  $B$ and $T$ in $P_{1+2k, 4+2k}$.  Since $\pi_1(T)\to \pi_1(P_{1+2k, 4+2k})=\ZZ$ is surjective,   $\pi_1(P_{1+2k, 4+2k}-T)\to
  \pi_1(P_{1+2k, 4+2k})$ is an isomorphism,  and $\pi_1(T_1)\to \pi_1(B)$ has image a cyclic summand, the gluing map for the symplectic sum can be chosen so that $B-nbd(T_1)\subset Y$ induces an isomorphism on fundamental groups. Hence $\pi_1(T_2)\to \pi_1(Y)$ is an isomorphism.

The symplectic sum
$$X_{93+2k, 94+2k}=Y\#_{T_2, T_4}S$$
  is a simply connected minimal symplectic 4-manifold with $e=189+4k$ and $\sigma=-1$, for any $k\ge 2$. \end{proof}

 \medskip

 Since any  symplectic signature zero 4-manifold has $e$ a multiple of $4$, there remain $45$ signature zero minimal symplectic 4-manifolds with odd intersection form to be constructed. Also missing are 48 signature $-1$ minimal symplectic 4-manifolds. Hence to complete the geography problem for minimal simply connected symplectic 4-manifolds of non-negative signature  there remain 97 manifolds to discover.

%**************************************************

%******* End of document **************************

%**************************************************

\end{document}